\documentclass[11pt]{amsart}

\usepackage{amsmath,amssymb}
\usepackage{fullpage}

\def\eps{\varepsilon}

\def\COMMENT#1{}
\def\TASK#1{}

\def\noproof{{\unskip\nobreak\hfill\penalty50\hskip2em\hbox{}\nobreak\hfill%
        $\square$\parfillskip=0pt\finalhyphendemerits=0\par}\goodbreak}
\def\endproof{\noproof\bigskip}
\newdimen\margin   
\def\textno#1&#2\par{%
    \margin=\hsize
    \advance\margin by -4\parindent
           \setbox1=\hbox{\sl#1}%
    \ifdim\wd1 < \margin
       $$\box1\eqno#2$$%
    \else
       \bigbreak
       \hbox to \hsize{\indent$\vcenter{\advance\hsize by -3\parindent
       \sl\noindent#1}\hfil#2$}%
       \bigbreak
    \fi}

\newtheorem{firstthm}{Proposition}
\newtheorem{thm}[firstthm]{Theorem}
\newtheorem{prop}[firstthm]{Proposition}
\newtheorem{lemma}[firstthm]{Lemma}

\newtheorem{defin}[firstthm]{Definition}
\newtheorem{conj}[firstthm]{Conjecture}
\newtheorem{claim}[firstthm]{Claim}
\newtheorem{fact}[firstthm]{Fact}

\begin{document}
\title{Embedding spanning bipartite graphs of small bandwidth}
\author{Fiachra Knox and Andrew Treglown}
\date{\today} 
\begin{abstract}  B\"ottcher, Schacht and Taraz~\cite{bot} gave a  condition on the minimum degree of a graph $G$ on 
$n$ vertices that ensures $G$ contains  every $r$-chromatic graph $H$ on $n$ vertices of bounded degree and of bandwidth $o(n)$, thereby
proving a conjecture of Bollob\'as and Koml\'os~\cite{komlos}. We strengthen this result in the case when $H$ is bipartite.
Indeed, we give an essentially best-possible condition on the degree sequence of a graph $G$ on $n$ vertices
that forces $G$ to contain every bipartite graph $H$ on $n$ vertices of bounded degree and of bandwidth $o(n)$. This also implies
an Ore-type result.
In fact, we prove a much stronger result where the condition on $G$ is relaxed to a certain robust expansion property.
Our result also confirms the bipartite case of a conjecture of Balogh, Kostochka and Treglown~\cite{bkt}
concerning the degree sequence of a graph which forces a perfect $H$-packing.
\end{abstract}
 
\maketitle

\section{Introduction}
A central problem in graph theory is to establish conditions on a graph $G$ which ensure that $G$ contains another graph $H$ as a 
spanning subgraph. Perhaps the best-known example of such a problem is when $H$ is a Hamilton cycle. Dirac's theorem~\cite{dirac}
states that any graph $G$ on $n$ vertices with minimum degree $\delta (G) \geq n/2$ contains a Hamilton cycle. 
The P\'osa--Seymour conjecture (see~\cite{posa} and~\cite{seymour}) states that any graph $G$ on $n$ vertices with $\delta (G) \geq
rn/(r+1)$ contains the $r$th power of a Hamilton cycle. (The $r$th power of a Hamilton cycle $C$ is obtained from $C$ by adding an edge
between every pair of vertices of distance at most $r$ on $C$.)  Koml\'os, S\'ark\"ozy and Szemer\'edi~\cite{kss} proved this conjecture
 for sufficiently large graphs.

There has also been significant attention on establishing minimum degree conditions which ensure a graph contains a
\emph{perfect $H$-packing}: Given a graph $H$, a perfect $H$-packing
in a graph $G$ is a collection of vertex-disjoint copies of $H$ which covers all the
vertices in $G$. (Perfect $H$-packings are also referred to as \emph{$H$-factors} or \emph{perfect $H$-tilings}.)
A seminal result in the area is the Hajnal--Szemer\'edi theorem~\cite{hs} which states that every graph $G$ whose order $n$ is
divisible by $r$ contains a perfect $K_r$-packing provided that $\delta (G) \geq (r-1)n/r$. (Corr\'adi and Hajnal~\cite{corradi}
had earlier
proved this result in the case when $r=3$.) 
Notice that in the case when $r+1$ divides $|G|$, a necessary condition for a graph $G$ to contain the $r$th
power of a Hamilton cycle is that $G$ contains a perfect $K_{r+1}$-packing. Thus, the P\'osa--Seymour conjecture implies
the Hajnal--Szemer\'edi theorem.
K\"uhn and Osthus~\cite{kuhn, kuhn2} characterised, up to an additive constant, the minimum degree which ensures a graph $G$ 
contains a perfect $H$-packing for an arbitrary graph $H$. (This improved previous bounds of Alon and Yuster~\cite{alon}
and  Koml\'os,  S\'ark\"ozy and Szemer\'edi~\cite{kss3}.)

It is desirable to find conditions that ensure a graph $G$ contains $H$ as a spanning subgraph where $H$ is \emph{any} 
graph from a large collection of graphs. That is, rather than finding individual results for specific graphs $H$, one seeks
more general, wide-reaching results. A graph $H$ on $n$ vertices is said to have \emph{bandwidth at most $b$}, if there exists a 
labelling of the vertices of $H$ by the numbers $1, \dots ,n$ such that for every edge $ij\in E(H)$ we have $|i-j|\leq b$.
Clearly every graph $H$ has bandwidth at most $|H|-1$. Thus, a perfect $H$-packing has bandwidth at most $|H|-1$.
Further, a Hamilton cycle has bandwidth $2$, and in general the $r$th power of a Hamilton cycle has bandwidth at most $2r$.
\COMMENT{If $|H|\leq 2r+1$ then bandwidth is $|H|-1$. Otherwise bandwidth at least $2r$: Vertex given label $1$ has $2r$ neighbours,
so one neighbour given label at least $2r+1$.
Let $H = 12 \ldots k1$ be a Hamilton cycle, and let $H^r$ be its $r$th power.  
Provided $|H|< 2r+1$ we are fine.
Otherwise, if we order the vertices $1, 2, k-1, 3, k-2, \ldots$, then for any edge $xy \in H^r$, the distance between $x$ and $y$ in this ordering is at most $2r$.}
B\"ottcher, Preussmann,  Taraz and  W\"urfl~\cite{BPTW10} proved that every planar graph  $H$ on $n$ vertices with bounded maximum degree has bandwidth at most $O(n/\log n)$.

The following result of B\"ottcher, Schacht and Taraz~\cite{bot} gives a condition on the minimum  degree of a graph $G$ on 
$n$ vertices
that ensures $G$ contains \emph{every} $r$-chromatic graph on $n$ vertices of bounded degree and of bandwidth $o(n)$, thereby
proving a conjecture of Bollob\'as and Koml\'os~\cite{komlos}.
\begin{thm}[B\"ottcher, Schacht and Taraz~\cite{bot}]\label{bst}
Given any $r,\Delta \in \mathbb N$ and any $\gamma >0$, there exist constants $\beta >0$ and $n_0 \in \mathbb N$
such that the following holds. Suppose that $H$ is an $r$-chromatic graph on $n \geq n_0$ vertices with $\Delta (H) \leq \Delta$ 
and bandwidth at most $\beta n$. If $G$ is a graph on $n$ vertices with
$$\delta (G) \geq \left( \frac{r-1}{r}+\gamma \right)n,$$
then $G$ contains a copy of $H$. \endproof
\end{thm}
Prior to the proof of Theorem~\ref{bst}, Csaba~\cite{csaba} and H\`an~\cite{han} proved the case when $H$ is bipartite and B\"ottcher, Schacht and 
Taraz~\cite{BST08} proved the case when $\chi(H)=3$. In this paper our focus is on strengthening Theorem~\ref{bst} in the case when $H$ is bipartite.

\subsection{Degree sequence conditions}
Dirac's theorem and the Hajnal--Szemer\'edi theorem are best possible in the sense that the minimum degree conditions in both these
results cannot be lowered. However, this does not mean that one cannot strengthen these results. Indeed, Chv\'atal~\cite{ch}
gave a condition on the degree sequence of a graph which ensures Hamiltonicity:
Suppose that the degrees of the graph $G$ are $d_1 \leq \dots \leq d_n$. If $n \geq 3$ and $d_{i} \geq i+1$ 
or $d_{n-i} \geq n-i$ for all $i < n/2$ then $G$ is Hamiltonian. 
Notice that  Chv\'atal's theorem is much stronger than Dirac's theorem since it allows for almost half of the vertices of $G$ to have degree less than $n/2$. 

Balogh, Kostochka and Treglown~\cite{bkt}
proposed the following two conjectures concerning the degree sequence of a graph which forces a perfect $H$-packing.
\begin{conj}[Balogh, Kostochka and Treglown~\cite{bkt}]\label{conj1}
Let $n,r \in \mathbb N$ such that $r $ divides $n$. Suppose that $G$ is a graph on $n$ vertices with degree sequence $d_1 \leq 
\dots \leq d_n$ such that:
\begin{itemize}
\item $d_i \geq (r-2)n/r+i $ for all $i < n/r$;
\item $d_{n/r+1} \geq (r-1)n/r$.
\end{itemize}
Then $G$ contains a perfect $K_r$-packing.
\end{conj}
Note that Conjecture~\ref{conj1}, if true, is much stronger than the Hajnal--Szemer\'edi theorem since the
degree condition allows for $n/r$ vertices to have degree less than $(r-1)n/r$.
\begin{conj}[Balogh, Kostochka and Treglown~\cite{bkt}]\label{conj2} Suppose $\gamma >0$ and $H$ is a graph with $\chi (H)=r$. Then there exists an integer $n_0=n_0 (\gamma ,H)$
such that the following holds. If $G$ is a graph whose order $n \geq n_0$ is divisible by $|H|$, and whose degree
sequence $d_1\leq \dots \leq d_n$ satisfies
\begin{itemize}
\item  $d_i \geq (r-2)n/r+i +\gamma n  $ for all $i < n/r$,
\end{itemize}
then $G$ contains a perfect $H$-packing. 
\end{conj}
In this paper we prove the following result which gives a condition on the degree sequence of a graph $G$ on 
$n$ vertices
that ensures $G$ contains every bipartite graph on $n$ vertices of bounded degree and of bandwidth $o(n)$.
\begin{thm}\label{mainthm1}
Given any $\Delta \in \mathbb N$ and any $\gamma >0$, there exists constants $\beta >0$ and $n_0 \in \mathbb N$
such that the following holds. Suppose that $H$ is a bipartite graph on $n \geq n_0$ vertices with $\Delta (H) \leq \Delta$ 
and bandwidth at most $\beta n$. Let $G$ be a graph on $n$ vertices with degree sequence $d_1\leq \dots \leq d_n$. If
\begin{itemize}
\item $d_{i} \geq i+\gamma n \ \text{ or } \ d_{n-i-\gamma n} \geq n-i \ \text{ for all } \ i < n/2$
\end{itemize}
then $G$ contains a copy of $H$.
\end{thm}
The degree sequence condition in Theorem~\ref{mainthm1} is similar to that in Chv\'atal's theorem, except that now we have 
two error terms in the condition.
Notice that Theorem~\ref{mainthm1} is much stronger than the bipartite case of Theorem~\ref{bst}.
Furthermore, in the case when $r=2$, Conjecture~\ref{conj2} is implied by Theorem~\ref{mainthm1}.

Theorem~\ref{mainthm1} is, up to the error terms, best-possible for many graphs $H$. Indeed, suppose that 
$H$ is a bipartite graph on an even number $n$ of vertices that contains a perfect matching. 
Suppose that $m \in \mathbb N$ such that $m<n/2$. Let
$G$ be a graph on $n$ vertices with vertex classes $V_1,V_2,V_3$ of sizes $m$, $m-1$ and $n-2m+1$ respectively and 
whose edge set contains all possible edges except for those in $V_1$ and between $V_1$ and $V_3$.
 Let $d_1 \leq \dots \leq d_n$ denote the
degree sequence of $G$. Then
\begin{itemize}
\item $d_{i} \geq i-1 \ \text{ and } \ d_{n-i+2} \geq n-i \ \text{ for all } \ i < n/2$,
\end{itemize}
but since $|V_1|>|V_2|$, $G$ does not contain a perfect matching and therefore $H$.

\subsection{Ore-type degree conditions}
Ore-type degree conditions consider the sum of the degrees of non-adjacent vertices of a graph. The name comes from Ore's theorem~\cite{ore}, which states that a graph~$G$ of
order $n\ge 3$ contains a
Hamilton cycle if $d(x)+d(y) \geq n$ for all non-adjacent $x \not = y\in V(G)$.
Recently, Ch\^au~\cite{chau} proved an Ore-type analogue of the P\'osa--Seymour conjecture in the case of the square of a Hamilton cycle (i.e. when $r=2$).

The following Ore-type result of Kierstead and Kostochka~\cite{kier} implies the Hajnal--Szemer\'edi theorem:
Let $n,r \in \mathbb N$ such that $r$ divides $n$. Suppose that
$G$ is a graph on $n$ vertices such that for all non-adjacent $x \not = y \in V(G)$,
$d(x)+d(y) \geq 2(r-1)n/r-1$.
Then $G$ contains a perfect $K_r$-packing.
K\"uhn, Osthus and Treglown~\cite{kotore}  characterised, asymptotically, the Ore-type degree condition which ensures a graph $G$ 
contains a perfect $H$-packing for an arbitrary graph $H$.

It is natural to seek an Ore-type analogue of Theorem~\ref{bst}. The following result provides such an analogue in the case when
$H$ is bipartite.
\begin{thm}\label{mainthm2}
Given any $\Delta \in \mathbb N$ and any $\gamma >0$, there exists constants $\beta >0$ and $n_0 \in \mathbb N$
such that the following holds. Suppose that $H$ is a bipartite graph on $n \geq n_0$ vertices with $\Delta (H) \leq \Delta$ 
and bandwidth at most $\beta n$. Let $G$ be a graph on $n$ vertices such that, for all non-adjacent $x \not = y \in V(G)$,
$$d(x)+d(y) \geq (1+\gamma ) n.$$
Then $G$ contains a copy of $H$.
\end{thm}
In Section~\ref{implysec} we show that Theorem~\ref{mainthm2} is a direct consequence of Theorem~\ref{mainthm1}.
Note that Theorem~\ref{mainthm2} is best-possible up to the error term for bipartite graphs $H$ on $n$ vertices
which do not contain an
isolated vertex. Indeed, let $G$ consist of a copy of $K_{n-1}$ and an isolated vertex. Then $G$ does not contain
$H$ but
$d(x)+d(y)=n-2$ for all non-adjacent $x \not = y \in V(G)$.

In light of Theorem~\ref{mainthm2}, we propose the following Ore-type analogue of Theorem~\ref{bst}.
\begin{conj}\label{conjKT}
Given any $r,\Delta \in \mathbb N$ and any $\gamma >0$, there exists constants $\beta >0$ and $n_0 \in \mathbb N$
such that the following holds. Suppose that $H$ is an $r$-chromatic graph on $n \geq n_0$ vertices with $\Delta (H) \leq \Delta$ 
and bandwidth at most $\beta n$. Let  $G$ be a graph on $n$ vertices such that, for all non-adjacent $x \not = y \in V(G)$,
$$d(x)+d(y)  \geq 2\left( \frac{r-1}{r}+\gamma \right)n.$$
Then $G$ contains a copy of $H$.
\end{conj}
If true, Conjecture~\ref{conjKT} is stronger than Theorem~\ref{bst}. B\"ottcher and M\"uller~\cite{botphd, bom} have proved
 the conjecture in the case when $r=3$.

\subsection{Robustly expanding graphs}\label{intro3}
An important and well-studied notion in graph theory is \emph{graph expansion}. We will consider the following
stronger notion of `robust expansion'. Roughly speaking, a graph $G$ on $n$ vertices is a robust expander if,
for every `reasonably sized' set $S \subseteq V(G)$, $G$ contains  at least $|S|+o(n)$ vertices that are adjacent to `many' vertices
in $S$. More formally,
let $0 < \nu \leq \tau < 1$. Suppose that $G$ is a graph on $n$ vertices and $S \subseteq V(G)$. Then the \emph{$\nu$-robust neighbourhood $RN_{\nu, G}(S)$ of $S$} is the set of vertices $v \in V(G)$ such that $|N(v) \cap S| \geq \nu n$.
We say that $G$ is a \emph{robust} $(\nu, \tau)$-\emph{expander} if every $S \subseteq V(G)$ with $\tau n \leq |S| \leq (1-\tau)n$ satisfies $|RN_{\nu, G}(S)| \geq |S| + \nu n$. 

The notion of robustly expanding (di)graphs was first introduced by K\"uhn, Osthus and Treglown in~\cite{KOT10}. 
\COMMENT{Shall we mention the applications this notion has already had?}
The following
result is an immediate consequence of Theorem~16 from~\cite{KOT10}.
\begin{thm}[K\"uhn, Osthus and Treglown~\cite{KOT10}]\label{expanderthm}
Given positive constants $\nu\le \tau\ll\eta<1$ there exists a positive integer~$n_0$
such that the following holds. Let~$G$ be a graph on~$n\ge n_0$ vertices with
$\delta(G)\ge \eta n$ which is a robust $(\nu,\tau)$-expander. Then~$G$ contains a Hamilton cycle. \endproof
\end{thm}
(Throughout the paper, we write $0<\alpha \ll \beta \ll \gamma$ to mean that we can choose the constants
$\alpha, \beta, \gamma$ from right to left. More
precisely, there are increasing functions $f$ and $g$ such that, given
$\gamma$, whenever we choose some $\beta \leq f(\gamma)$ and $\alpha \leq g(\beta)$, all
calculations needed in our proof are valid. 
Hierarchies of other lengths are defined in the obvious way.)

We will use Theorem~\ref{expanderthm} to prove the following result concerning embedding bipartite graphs of small bandwidth.
\begin{thm}\label{mainthm}
Given $\Delta \in \mathbb N$ and positive constants $\nu\le \tau\ll\eta<1$ there exist constants $\beta >0$ and~$n_0 \in \mathbb N$
such that the following holds.
Suppose that $H$ is a bipartite graph on $n \geq n_0$ vertices with $\Delta (H) \leq \Delta$ 
and bandwidth at most $\beta n$. Let $G$ be a graph on $n$ vertices  with
$\delta(G)\ge \eta n$ which is a robust $(\nu,\tau)$-expander.
Then $G$ contains a copy of $H$.
\end{thm}
In Section~\ref{implysec} we show that Theorem~\ref{mainthm} implies Theorem~\ref{mainthm1} and that Theorem~\ref{mainthm1}
implies Theorem~\ref{mainthm2}. Thus, we only
prove Theorem~\ref{mainthm} directly.

Note that Theorem~\ref{mainthm} is very general in the sense that it allows for the graph $G$ to have small minimum degree
(although $\delta (G)$ must be linear). Furthermore, there are examples of graphs $G$ that satisfy the hypothesis of
Theorem~\ref{mainthm} and whose \emph{maximum degree} is  also small. Indeed, let $0<\nu \ll \tau \ll \eta < 1$ such that
$1/\eta$ is an odd integer. Further choose $n \in \mathbb N$ such that $\eta n \in \mathbb N$. 
Define $G$ to be the blow-up of a cycle on $1/\eta$ vertices, such that each vertex class of $G$ contains $\eta n$ vertices.
Thus, $|G|=n$ and $\delta (G)=\Delta (G)=2\eta n$. It is easy to check that $G$ is a robust $(\nu,\tau)$-expander.
\COMMENT{For any set $S$ such that $\tau n \leq |S| \leq (1 - \tau)n$, let $k$ be such that $ (k-1)\eta  n+(\nu n)/\eta
 < |S| \leq  k\eta n+(\nu n)/ \eta$. Then there must be at least $k$ vertex classes $U$ such that $|S \cap U| \geq \nu n$, and for each such $U$, both of the neighbouring clusters of $U$ must be contained in $RN_{\nu, G}(S)$. 
 If $k=1/\eta$ then we are done. Otherwise, by a careful induction argument one can show that
there must be at least $k+1$ such `neighbour' clusters. Hence $|RN_{\nu, G}(S)| \geq (k+1)\eta n \geq k\eta n + (\nu n)/\eta +
\nu n \geq |S| + \nu n$.}
Given constants $0<\nu \ll \tau \ll p <1$, with high probability $G(n,p)$ is a robust $(\nu,\tau)$-expander with minimum degree
at least $pn/2$ and maximum degree at most $2pn$.

Theorem~\ref{mainthm} therefore implies that, with high probability, $G(n,p)$ contains all  bipartite graphs $H$ on $n$ vertices
of bounded degree and bandwidth $o(n)$. A result of Huang, Lee and Sudakov~\cite{huang} actually implies that, with high
probability, any spanning subgraph $G'$ of $G(n,p)$ with minimum degree $\delta (G') \geq (1/2+o(1))np$  contains all such $H$.

\section{Notation and preliminaries}
\subsection{Notation}
Throughout this paper we omit floors and ceilings whenever this does not affect the
argument.
We write  $|G|$ for the  order of a graph $G$,
$\delta (G)$
and $\Delta (G)$ for its minimum and maximum degrees respectively and $\chi (G)$ for
its chromatic number. The degree of a vertex $x \in V(G)$ is denoted by $d(x)$ and its neighbourhood
by $N(x)$. Given $S \subseteq V(G)$ we define $N(S):=\bigcup_{v \in S} N(v)$.

Given disjoint $A,B \subseteq V(G)$ the number of  edges with one endpoint in $A$ and one endpoint in $B$ 
is denoted by $e_G
(A,B)$. We write $(A,B)_G$ for the bipartite subgraph of~$G$
with vertex classes $A$ and $B$ whose edges are precisely those edges in $G$ with one endpoint in $A$ and the other in $B$.
Often we will write $(A,B)$, for example, if this is unambiguous.

\subsection{Degree sequence and Ore-type conditions forcing robust expansion}\label{implysec}
The following result is an immediate consequence of Lemma~13 from~\cite{KOT10}.
\begin{lemma}[\cite{KOT10}]\label{lemimply}
Given positive constants $\tau \ll \eta < 1$ there exists an integer~$n_0$ such that whenever~$G$
is a graph on $n\ge n_0$ vertices with
\begin{align*}
 d _i\geq i + \eta n  \text{ or } d _{n-i-\eta n} \geq n-i \text{ for all }i <n/2,
\end{align*}
 then $\delta (G) \ge \eta n$ and~$G$ is a robust $(\tau^2,\tau)$-expander. \endproof
\end{lemma}
Notice that Lemma~\ref{lemimply} together with Theorem~\ref{mainthm} implies Theorem~\ref{mainthm1}.
We now show that Theorem~\ref{mainthm1} implies Theorem~\ref{mainthm2}.

\begin{lemma} Let $ \gamma >0 $. Suppose $G$ is a graph on $n$ vertices such that, for all non-adjacent $x \not = y \in V(G)$,
$$d(x)+d(y) \geq (1+2\gamma ) n.$$ Let $d_1\leq \dots \leq d_n$ denote the degree sequence of
$G$. Then 
\begin{align*}
 d _i\geq i + \gamma n  \text{ or } d _{n-i-\gamma n} \geq n-i \text{ for all }i <n/2.
\end{align*}
\end{lemma}
\proof Firstly note that for $(1 - \gamma)n/2 \leq i < n/2 $ we wish to show that either $d_{n - i' - \gamma n} \geq n  - i'$ or $d_{i'} \geq i' + \gamma n$, where $i' := n - i - \gamma n$. Notice that $n/2 -\gamma n < i' \leq n/2 -\gamma n/2$. Thus, 
it suffices to only consider $i$ such that $1 \leq i \leq (1 - \gamma)n/2$.

Suppose there is some $1\leq i \leq (1-\gamma)n/2$ such that the statement does not hold. Then there is a set $A$ of $i $ vertices, each of degree less than $i + \gamma n \leq n/2 + \gamma n/2$. So for any $x, y \in A$, $d(x) + d(y) <(1 + 2\gamma)n$ and hence $G[A]$ is a clique. Set $B := V(G) \backslash A$. Note that $e_G(A, B) < (\gamma n +1)i$. 
Hence, there is a vertex $x \in B$ that receives less than $\min \{ \gamma n +1 , i \}$ edges from $A$. Therefore, there is a vertex
$y \in A$ such that $xy \not \in E(G)$.  Thus, $d(x) + d(y) < (n-i-1+ \gamma n+1) +(i+\gamma n)\leq (1+2 \gamma)n$, contradicting our assumption. 
\endproof

\section{Outline of the proof of Theorem~\ref{mainthm}} \label{SecSketch}
\subsection{Proof overview}
The overall strategy is similar to that of the proof of Theorem~\ref{bst} in~\cite{bot}. Indeed, as in~\cite{bot} the proof
is split into two main lemmas; the \emph{Lemma for $G$} and the \emph{Lemma for $H$}. However, many of the methods used in~\cite{bot}
break down in our setting so our argument proceeds somewhat differently. 

The role of the Lemma for $G$ (Lemma~\ref{LemmaForG}) is to obtain some special structure within $G$ so that it will be suitable for
embedding $H$ into; By applying Theorem~\ref{expanderthm}, 
we show that $G$ contains a spanning subgraph $G'$ which `looks' like the blow-up of a cycle 
$C=V_1V_2\dots V_{2k}V_1$.
More precisely, there is a partition $V_1, \dots , V_{2k}$ of $V(G)$ such that:
\begin{itemize}
\item[(i)] $(V_{2i-1},V_{2i})_{G'}$ is a `super-regular' pair of density at least $d>0$ for each $1\leq i \leq k$;
\item[(ii)] $(V_{2i},V_{2i+1})_{G'}$ is an `$\eps$-regular' pair of density at least $d$ for each $1\leq i \leq k$.

\end{itemize}
Furthermore, there are  even integers $1\leq i_1 \not = j_1 \leq 2k$
such that:
\begin{itemize}
\item[(iii)]  $(V_{i_1},V_{j_1})_{G'}$ is `$\eps$-regular' with density at least $d$.
\end{itemize}
(So $V_{i_1}V_{j_1}$  can be thought of as a chord of $C$.)
Crucially, this partition is `robust' in the sense that one can modify the sizes of each partition class $V_i$ somewhat without
destroying the  properties (i)--(iii). (This is made precise by the \emph{Mobility lemma} given in Section~\ref{SecMobility}.)

Set $c:=V_{i_1}V_{j_1}$. The role of the Lemma for $H$ (Lemma~\ref{LemmaForH}) is to construct a graph homomorphism $f$ from $H$
to $C \cup \{c\}$ in such a way that `most' of the edges of $H$ are mapped to edges of the form $V_{2i-1}V_{2i}$ for some $i$.
(Recall that these are the edges which correspond to super-regular pairs in $G'$.)
The homomorphism $f$ is such that every $V_i \in C$ receives roughly $|V_i|$ vertices of $H$. So $f$ can be viewed as a `guide'
as to which vertex class $V_i \subseteq V(G)$ each vertex from $H$ is embedded into. In particular, since the partition
$V_1, \dots , V_{2k}$ is `robust', we can alter the sizes of the classes $V_i$ such that (i)--(iii) still hold and  so that
now $|f^{-1} (V_i)|=|V_i|$ for all $i$.
Properties (i)--(iii) then allow us to apply the Blow-up lemma~\cite{kss2} to embed $H$ into $G'$ and thus $G$.
(Actually we apply a result from~\cite{botphd} which is a consequence of the Blow-up lemma.)

\subsection{Techniques for the Lemma for $G$} In order to obtain the partition $V_1, \dots , V_{2k}$ 
of $V(G)$ we modify a partition $V'_0, V'_1, \dots , V'_{2k}$
obtained by applying Szemer\'edi's Regularity
lemma~\cite{reglem} to $G$. Roughly speaking, $V'_1, \dots ,V'_{2k}$ will satisfy (i)--(iii). Thus, we need to redistribute 
the vertices from $V'_0$ into the other vertex classes whilst retaining these properties. 
We also require our partition $V_1, \dots, V_{2k}$ to satisfy
$|V_{2i-1}| \approx |V_{2i}|$ for each $1\leq i \leq k$. So we need to redistribute vertices in a `balanced' way.
In the Lemma for $G$ in~\cite{bot}, the minimum degree condition of Theorem~\ref{bst} is heavily relied on
 to achieve this. However, our graph
$G$ may have very small minimum degree. So instead we introduced the notion of a `shifted $M$-walk' to help us redistribute vertices:
Given a perfect matching $M$ in a graph $R$ a shifted $M$-walk is a walk whose edges alternate between edges of $M$ 
and edges of $R\backslash M$ (see Section~\ref{SecShiftedWalks} for the precise definition). Since $G$ is a robust expander,
we can find short shifted $M$-walks in a reduced graph $R$ of $G$. 
(Here, $M$ will be the perfect matching in $R$ that corresponds to the super-regular pairs from (i) above.)
These walks act as a `guide' as to how we redistribute vertices
amongst the vertex classes. 

\subsection{Techniques for the Lemma for $H$}
In~\cite{bot} the techniques used are actually strong enough to prove a more general result than Theorem~\ref{bst} (and so
Theorem~\ref{bst} is not proved directly). For example,
in the case when $r=2$, their result concerns not only bipartite $H$ but also a special class of $3$-colourable graphs $H$ where
the third colour class is very small (see Theorem~2 in~\cite{bot} for precise details). One example of such a graph $H$
is a Hamilton cycle $C'$ with a chord between two vertices of distance $2$ on $C'$. $H$ is $3$-colourable and has bounded
bandwidth. However, $H$ cannot be embedded into \emph{every} graph $G$ satisfying the hypothesis of Theorem~\ref{mainthm}.
Indeed, consider the graph $G$ defined at the end of Section~\ref{intro3}. 

In particular, this means we have to approach the proof of the Lemma for $H$ differently: Since $H$ has bandwidth $o(n)$
we can chop $V(H)$ into small linear sized segments $A_1, B_1, \dots , A_m ,B_m$ where all the edges of $H$ lie in pairs of
the form $(A_i,B_i)_H$ and $(B_i, A_{i+1})_H$ and such that $A:=\cup ^m _{i=1} A_i$ and $B:=\cup ^m _{i=1} B_i$ are the colour
classes of $H$. Ideally we would want to construct $f$ to map the vertices of $A_1$ into $V_1$, the vertices of $B_1$ into $V_2$
and so on, continuing around $C$ many times until all the vertices have been assigned. However, since $|A|$ and $|B|$ may vary
widely, this would map vertices in an unbalanced way. That is, the total number of vertices mapped to `odd' classes $V_{2i-1}$ 
would  differ widely from the total number of vertices mapped to `even' classes $V_{2i}$. We get around this problem by
using the chord $c=V_{i_1}V_{j_1}$ to `flip' halfway in the process. So after this, vertices from the $B_i$ are mapped to
`odd' classes  $V_{2i-1}$  and vertices from the $A_i$ are mapped to the `even' classes $V_{2i}$. We also `randomise' part of the
mapping procedure to ensure that the number of vertices of $H$ assigned to each $V_i$ is approximately $|V_i|$.
(A randomisation technique of a similar flavour was used in~\cite{KMOprep}.)

\section{The Regularity lemma } \label{SecRegularity}
In the proof of the Lemma for $G$ (Lemma~\ref{LemmaForG}) we will use Szemer\'edi's Regularity
lemma~\cite{reglem}.
In this section we will introduce all the information we require about this
result.
To do this, we firstly introduce some more notation.
The \emph{density} of a bipartite graph $G$ with vertex classes~$A$ and~$B$ is
defined to be
$$d(A,B):=\frac{e(A,B)}{|A||B|}.$$
Given any $\eps, d>0$, we say that $G$ is \emph{$(\eps,d)$-regular} if $d(A,B)\geq d$ and, for all sets
$X \subseteq A$ and $Y \subseteq B$ with $|X|\ge \eps |A|$ and
$|Y|\ge \eps |B|$, we have $|d(A,B)-d(X,Y)|< \eps$. 
We say that $G$ is \emph{$(\eps,d)$-super-regular} if additionally every vertex $a \in A$ has at least $d|B|$ neighbours
in $B$ and every vertex $b \in B$ has at least $d|A|$ neighbours in $A$.
We also say that $(A,B)$ is an \emph{$(\eps,d)$-(super-)regular pair}.
We will frequently use the following simple fact.

\begin{fact} \label{RegularDegree} Let $\eps , d >0$. Suppose that $G = (A, B)$ is an $(\varepsilon, d)$-regular pair. Let $B' \subseteq B$ be such that $|B'| \geq \varepsilon |B|$. Then there are at most $\varepsilon |A|$ vertices in $A$ with fewer than $(d - \varepsilon)|B'|$ neighbours in $B'$. \endproof \end{fact}

We  will also require the next simple proposition which allows us to modify a (super-)regular pair without destroying its (super-)regularity (see e.g., \cite[Proposition 8]{BST08}).
\begin{prop} \label{PreserveReg} Let $(A, B)$ be an $(\varepsilon, d)$-regular pair, and let $A'$ and $B'$ be vertex sets with $|A' \Delta A| \leq \alpha |A|$ and $|B' \Delta B| \leq \beta |B|$. Then $(A', B')$ is an $(\varepsilon ', d')$-regular pair where
$$\varepsilon ' := \varepsilon + 3(\sqrt{\alpha} + \sqrt{\beta}) \text{ and } d' := d - 2(\alpha + \beta).$$
If, moreover, $(A, B)$ is $(\varepsilon, d)$-super-regular and each vertex in $A'$ has at least $d|B'|$ neighbours in $B'$ and each vertex in $B'$ has at least $d|A'|$ neighbours in $A'$, then $(A', B')$ is $(\varepsilon', d')$-super-regular. \endproof
\end{prop}

We will use the following degree form of Szemer\'edi's Regularity lemma~\cite{reglem} which can be
easily derived from the classical version.
\begin{lemma}[Regularity lemma] \label{Regularity} For every $\varepsilon > 0$ and $k_0 \in \mathbb{N}$ there exists $K_0 = K_0(\varepsilon, k_0)$ such that for every $d \in [0, 1]$ and for every graph $G$ on $n \geq K_0$ vertices there exists a partition $V_0, V_1, \ldots, V_k$ of $V(G)$ and a spanning subgraph $G'$ of $G$, such that the following conditions hold:
\begin{itemize}
\item [\rm (i)] $k_0 \leq k \leq K_0$,
\item [\rm (ii)] $d_{G'}(x) \geq d_G(x) - (d + \varepsilon) n$ for every $x \in V(G)$,
\item [\rm (iii)] the subgraph $G'[V_i]$ is empty for all $1 \leq i \leq k$,
\item [\rm (iv)] $|V_0| \leq \varepsilon n$,
\item [\rm (v)] $|V_1| = |V_2| = \ldots = |V_k|$,
\item [\rm (vi)] for all $1 \leq i < j \leq k$  either $(V_i, V_j)_{G'}$ 
is an $(\varepsilon, d)$-regular pair or $G'[V_i, V_j]$ is empty.
\end{itemize} \endproof
\end{lemma}
We call $V_1, \dots, V_k$ {\it clusters}, $V_0$ the {\it exceptional set} and the
vertices in~$V_0$ {\it exceptional vertices}. We refer to~$G'$ as the {\it pure graph}.
The {\it reduced graph~$R$ of~$G$ with parameters $\varepsilon$, $d$ and~$k_0$} is the graph whose 
vertices are $V_1, \dots , V_k$ and in which $V_i V_j$ is an edge precisely when $(V_i,V_j)_{G'}$
is $(\varepsilon,d)$-regular.

The following result implies that the property of a graph $G$ being a robust expander is `inherited' by the reduced graph $R$ of $G$.
 It is an 
immediate consequence of Lemma~14 from~\cite{KOT10}.
\begin{lemma}[\cite{KOT10}] \label{PreserveRE} Let $k_0,n_0$ be positive integers and let $\varepsilon,d,\eta,\nu,\tau$ be positive constants such that
$1/n_0\ll  \varepsilon \ll d\ll \nu,\tau,\eta < 1$ and such that $k_0\ll n_0$. Let~$G$ be a
graph on $n\ge n_0$ vertices with $\delta  (G) \ge \eta n$ and such that~$G$
is a robust $(\nu,\tau)$-expander. Let~$R$ be the
reduced graph of $G$ with parameters $\varepsilon$, $d$ and~$k_0$. Then $\delta (R) \ge \eta |R|/2$
and~$R$ is a robust $(\nu/2,2\tau)$-expander. \endproof
\end{lemma}

\section{Useful results}\label{useful}
\subsection{Shifted walks and robust expanders} \label{SecShiftedWalks}
\COMMENT{Some blah-blah about how these results are useful}
Let $G$ be a graph containing a perfect matching $M$. A \emph{shifted $M$-walk} in $G$ with endpoints $a = v_1$ and $b = v_{2\ell}$ is a walk $v_1 v_2 \ldots v_{2\ell}$ in $G$ 
such that $v_{2i} v_{2i+1} \in M$ for every $1 \leq i \leq \ell - 1$ and $v_{2i-1} v_{2i} \not \in M$ for any $1 \leq i \leq \ell $. A shifted $M$-walk is \emph{simple} if it contains each edge of $M$ at most twice.
Note that a path containing a single edge $v_1 v_2 \not \in M$ is a (simple) shifted $M$-walk for any perfect matching $M$.

\begin{lemma} \label{FindSimpleWalk}
Let $G$ be a graph containing a perfect matching $M$, and let $W$ be a shifted $M$-walk in $G$ with endpoints $a$ and $b$. Then $W$ contains a simple shifted $M$-walk $W'$ with endpoints $a$ and $b$.
\end{lemma}
\proof We proceed by induction. Let $W = v_1 \ldots v_{2 \ell}$. If $W$ is already simple then we set $W' := W$; otherwise, there exists an edge $xy \in M$ which appears at least three times in $W$. Let the first three appearances of $xy$ be $v_{i_1} v_{i_1 + 1}$, $v_{i_2} v_{i_2 + 1}$ and $v_{i_3} v_{i_3 + 1}$ (in order). Now each of $v_{i_1}$, $v_{i_2}$ and $v_{i_3}$ is either $x$ or $y$, and so without loss of generality we can assume that $v_{i_1} = v_{i_2} = x$. Now we have a shorter shifted $M$-walk $v_1 \ldots v_{i_1} v_{i_2 + 1} \ldots v_{2\ell}$. \endproof

\begin{lemma} \label{FindPureWalk} Let $M$ be a perfect matching in a graph $G$ 
and let $A \subseteq V(M)$ be a set containing at most one vertex from each edge of $M$. Suppose that
 $W$ is a shifted $M$-walk both of whose endpoints lie in $A$. Then $W$ contains a shifted $M$-walk $W'$ such that the endpoints of $W'$ both lie in $A$ and no other vertices of $W'$ lie in $A$.
\end{lemma}
\proof Let $W = v_1 \ldots v_{2\ell}$. Let $1 \leq i_1 \leq \ell$ be minimal such that $v_{2i_1} \in A$, and let $1 \leq i_2 \leq i_1$ be maximal such that $v_{2i_2 - 1} \in A$. Now $v_{2i_2 - 1} v_{2i_2} \ldots v_{2i_1 - 1} v_{2i_1}$ is the desired shifted $M$-walk. \endproof

In the proof of the Lemma for $G$ we will use shifted walks in the reduced graph $R$ of $G$ as a ``guide'' as to how 
to redistribute vertices in $G$. Since the reduced graph $R$ will be a robust expander, the following result ensures we can find
our desired shifted walks.

 Let $G$ be a graph containing a perfect matching $M$, and let $A \subseteq V(G)$. For each $v \in V(G)$, let $v' \in V(G)$ be the unique vertex such that $vv' \in M$. The \emph{shifted $M$-neighbourhood} of $A$ is the set $SN_M(A) = \{v' \mid v \in N(A)\}$. $SN_M^r(A)$ is defined recursively by $SN_M^1(A) := SN_M(A)$ and $SN_M^r(A) := SN_M(SN_M^{r-1}(A))$ for $r \geq 2$.
\begin{lemma} \label{FindWalk} Let $0 < \nu \leq \tau < \eta \ll 1$ be constants. Suppose $G$ is a graph on $n$ vertices 
with $\delta(G) \geq \eta n$ which is a robust $(\nu, \tau)$-expander, and let $M$ be a perfect matching in $G$. Then for any $a \in V(G)$, $G$ contains a shifted $M$-walk of length at most $3/\nu$ which both starts and finishes at $a$.
\end{lemma}
\proof The minimum degree condition implies that $|SN_M(a)| = |N(a)| \geq \eta n \geq \tau n$.  Since $G$ is a robust $(\nu, \tau)$-expander,
$$|SN_M^r(a)| = |N(SN_M^{r-1}(a))| \geq \min\left\{|SN_M^{r-1}(a)| + \nu n, (1 - \tau + \nu) n\right\},$$ 
 for all $r \geq 2$.
Hence $|SN_M^{1/2 \nu}(a)| \geq (\tau + 1/2 - \nu)n$ and so 
$$|N(SN_M^{1/2 \nu}(a))| \geq ( 1/2+\tau)n.$$ 
Thus, there exists some edge $vv' \in M$ such that both $v$ and $v'$ lie in $N(SN_M^{1/2 \nu}(a))$. This implies that there exists a shifted $M$-walk $P$ with endpoints $a$ and $v$ and a shifted $M$-walk $P'$ with endpoints $a$ and $v'$, each of length at most $1/\nu + 1$. Now $P \cup vv' \cup P'$ forms a shifted $M$-walk of length at most $2/\nu + 3 \leq 3/\nu$ which starts and finishes at $a$. 
\endproof

The next lemma allows us to delete a small number of vertices from a robust expander without destroying this property.

\begin{lemma} \label{PreserveRobustExp} Let $0 < \alpha < \nu \leq \tau \ll 1$ be constants. Suppose that $G$ is a graph
on $n$ vertices which is a robust $(\nu, \tau)$-expander and let $S \subseteq V(G)$ be a set of size $\alpha n$. Then $G - S$ is a robust $(\nu - \alpha, \tau + \alpha)$-expander.
\end{lemma}
\proof Let $G' := G-S$ and $n' := |G'|$. Consider any $A \subseteq V(G')$ such that $(\tau + \alpha)n' \leq |A| \leq (1 - \tau - \alpha)n'$. Set $A' := A \cup S$. Then 
$$\tau n \leq (\tau + \alpha)n' + \alpha n \leq |A'| \leq (1 - \tau  - \alpha)n' + \alpha n \leq (1 - \tau)n.$$ 
So $|RN_{\nu, G}(A')| \geq |A'| + \nu n$. Now every vertex of $RN_{\nu, G}(A')$ has at least $\nu n$ neighbours in $A'$ and since $|S| = \alpha n$, at least $(\nu - \alpha)n \geq (\nu - \alpha)n'$ of these must lie in $A$. Hence every vertex of $RN_{\nu, G}(A') \backslash S$ lies in $RN_{\nu - \alpha, G'}(A)$, and so $|RN_{\nu - \alpha, G'}(A)| \geq |A'| + \nu n - |S| \geq |A| + (\nu - \alpha) n'$, as desired. \endproof

\subsection{Probabilistic bounds}
The following two probabilistic bounds will be used in the proof of  the Lemma for $H$ (Lemma~\ref{LemmaForH}).
\begin{lemma}[\cite{KMOprep}, Lemma 2.1] \label{EvenBinomial} Suppose that $1/k \ll p, (1-p), \varepsilon$, that $n \geq k^3/6$, and that $X \sim Bin(n, p)$. Then for any $0 \leq r \leq k-1$,
$$\frac{1 - \varepsilon}{k} \leq \mathbb{P}(X \equiv r \mod k) \leq \frac{1 + \varepsilon}{k}.$$ \endproof
\end{lemma}

\begin{lemma}[\cite{SVprep}, Proposition 1.1] \label{Azuma} Let $X_1, \ldots, X_n$ be random variables taking values in $[0, 1]$, such that for each $1 \leq k \leq n$, 
$$\mathbb{E}[X_k \mid X_{k-1}, \ldots, X_1] \leq a_k.$$
Let $\mu := \sum_{k = 1}^n {a_k}$. Then for any $0 < \delta \leq 1$,
$$\mathbb{P}\left[\sum_{k=1}^n X_k > (1 + \delta)\mu\right] \leq e^{-\frac{\delta^2 \mu}{3}}.$$ \endproof
\end{lemma}
We also require the following expectation bound.
\begin{lemma} \label{ConditionalExp} Suppose that $X$ and $Y$ are integer-valued random variables and that $B$ is an event, such that for each $x, y \in \mathbb{Z}$, $\mathbb{P}[X = x \mid B \cap (Y = y)] = \mathbb{P}[X = x \mid Y = y]$. Then
$$\mathbb{E}[X \mid B] \leq \max_{y \in \mathbb{Z}} \mathbb{E}[X \mid Y = y].$$
\end{lemma}

\proof Note that for each $x \in \mathbb{Z}$,
$$\mathbb{P}[X = x \mid B] = \sum_{y \in \mathbb{Z}} \mathbb{P}[(X = x) \cap (Y = y) \mid B].$$
Further,
\begin{align*}
\mathbb{P}[(X = x) \cap (Y = y) \mid B] &= \frac{\mathbb{P}[(X = x) \cap (Y = y) \cap B]}{\mathbb{P}[B]} \\
&= \frac{\mathbb{P}[(X = x) \cap (Y = y) \cap B]}{\mathbb{P}[(Y = y) \cap B]} \cdot \frac{\mathbb{P}[(Y = y) \cap B]}{\mathbb{P}[B]} \\
&= \mathbb{P}[X = x \mid Y = y] \cdot \mathbb{P}[Y = y \mid B].
\end{align*}
Hence,
\begin{align*}
\mathbb{E}[X \mid B] &= \sum_{x \in \mathbb{Z}} x \mathbb{P}[X = x \mid B] = \sum_{x \in \mathbb{Z}} x \sum_{y \in \mathbb{Z}} \mathbb{P}[X = x \mid Y = y] \cdot \mathbb{P}[Y = y \mid B] \\
&= \sum_{y \in \mathbb{Z}} \mathbb{P}[Y = y \mid B] \sum_{x \in \mathbb{Z}} x \mathbb{P}[X = x \mid Y = y] = \sum_{y \in \mathbb{Z}} \mathbb{P}[Y = y \mid B] \cdot \mathbb{E}[X \mid Y = y] \\
&\leq \max_{y \in \mathbb{Z}} \mathbb{E}[X \mid Y = y].
\end{align*}
\endproof

\section{The Mobility lemma} \label{SecMobility}
In order to state our next result we first introduce a slight variant of the notion of a reduced graph.
Let $\eps, \eps ' , d, d'>0$. Suppose that $G$ is a graph and $V_1, \dots , V_k$ is a partition of $V(G)$. We say that a graph
$R$ is an \emph{$(\eps,d)$-reduced graph of $G$ on $V_1, \dots , V_k$} if the following holds:
\begin{itemize}
\item $V(R)=\{V_1, \dots , V_k \}$;
\item If $V_i V_j \in E(R)$ then $(V_i, V_j)_G$ is an $(\eps, d)$-regular pair (for all $1\leq i \not = j \leq k$).
\end{itemize}
Suppose $V'_1 , \dots , V'_k$ is another partition of $V(G)$
and $R$ is as above (in particular, $V(R)=\{V_1, \dots , V_k \}$). We also say that $R$ is an 
\emph{$(\eps ',d')$-reduced graph of $G$ on $V'_1, \dots , V'_k$} if the following holds:
\begin{itemize}
\item If $V_i V_j \in E(R)$ then $(V'_i, V'_j)_G$ is an $(\eps ', d')$-regular pair (for all $1\leq i \not = j \leq k$).
\end{itemize}

Suppose that $V_1 , \dots , V_k$ and $V'_1, \dots , V'_k$ are both partitions of the vertex set of a graph $G$. 
Given a cluster $V=V_i$ for some $1 \leq i \leq k$, we will often denote by $V'$ the cluster $V'_i$.

We will apply the next result in the proof of the Lemma for $G$ (Lemma~\ref{LemmaForG}) so that we can
 alter a particular partition of a graph $G$ somewhat without destroying the structure of our reduced graph $R$.

\begin{lemma}[Mobility lemma] \label{MobLemma} Let $k \in \mathbb{N}$, and let $\xi, \varepsilon, \varepsilon ', d', d$ be positive constants such that
$$0 < \xi \ll 1/k \ll \varepsilon \ll \varepsilon ' \ll d' \ll d \ll 1.$$
Suppose $G$ is a graph on $n$ vertices, $A_1, B_1, A_2, B_2, \ldots, A_k, B_k$ is a partition of $V(G)$ such that $|A_i|, |B_i| \geq n/3k$ for all $1 \leq i \leq k$ and $R$ is an $(\varepsilon, d)$-reduced graph on $A_1, B_1, \ldots, A_k, B_k$. Let $(a_i)_{i=1}^k$ and $(b_i)_{i=1}^k$ be integers. Suppose that the following conditions hold:
\begin{itemize}
\item[\rm (i)] $R$ contains the Hamilton cycle $C =  A_1 B_1 A_2 B_2 \ldots A_k B_k A_1$;
\item[\rm (ii)] $R$ contains an edge $A_{i_1} A_{j_1}$ for some $i_1 \neq j_1$;
\item[\rm (iii)] $R$ contains an edge $B_{i_2} B_{j_2}$ for some $i_2 \neq j_2$;
\item[\rm (iv)] The pair $(A_i, B_i)_G$ is $(\varepsilon, d)$-super-regular for all $1 \leq i \leq k$;
\item[\rm (v)] $|a_i|, |b_i| < \xi n$ for each $1 \leq i \leq k$;
\item[\rm (vi)] $\sum_{i=1}^k a_i + \sum_{i=1}^k b_i = 0$;
\item[\rm(vii)] $|\sum_{i=1}^k a_i| = |\sum_{i=1}^k b_i| \leq \xi n$.
\end{itemize}
Then there exists a partition $A'_1, B'_1, A'_2, B'_2, \ldots, A'_k, B'_k$ of $V(G)$ such that $|A'_i| = |A_i| + a_i$ and $|B'_i| = |B_i| + b_i$ for each $1 \leq i \leq k$, $R$ is an $(\varepsilon', d')$-reduced graph of $G$ on $A'_1, B'_1, A'_2, B'_2, \ldots, A'_k, B'_k$, and 
 $(A'_i, B'_i)_G$ is $(\varepsilon', d')$-super-regular for each $1 \leq i \leq k$.
\end{lemma}

\proof Without loss of generality we may assume that $\sum_{i=1}^k a_i \geq 0$. (As a consequence of this assumption we will in fact only need the edge $B_{i_2} B_{j_2}$, and not the edge $A_{i_1} A_{j_1}$.) Note that by (iii) and Fact \ref{RegularDegree} there are at least $(1 - \varepsilon)|B_{i_2}| \gg \xi n$ vertices in $B_{i_2}$ with at least $(d - \varepsilon)|B_{j_2}|$ neighbours in $B_{j_2}$. Pick $\sum_{i=1}^k a_i \leq \xi n$ of these vertices and move them from $B_{i_2}$ into $A_{j_2}$. Call the resulting sets $B^*_{i_2}$ and $A^*_{j_2}$ respectively.

We now perform an iterative procedure which will reassign vertices among the vertex classes $(A_i)_{i = 1}^k$ and, separately, $(B_i)_{i = 1}^k$. Initially we define the classes $A^*_i = A_i$ for each $i \neq j_2$ and $B^*_i = B_i$ for each $i \neq i_2$. Roughly speaking, $A^*_i$ (or $B^*_i$) will be the current version of $A_i$ (or $B_i$). The choice of how we defined $B^*_{i_2}$ and $A^*_{j_2}$ is such that,
initially,
\begin{align}\label{iter1}
\sum ^k _{i=1} |A^* _i|=\sum ^k _{i=1} |A_i|+ \sum_{i=1}^k a_i \ \text{ and } \
\sum ^k _{i=1} |B^* _i|=\sum ^k _{i=1} |B_i|+ \sum_{i=1}^k b_i.
\end{align}
Throughout the procedure we will ensure that (\ref{iter1}) holds. Furthermore, throughout we will ensure that
\begin{equation} \label{DifferenceBound} 
|A^*_i \Delta A_i|, |B^*_i \Delta B_i| \leq 5k\xi n \leq \eps  |A_i| , \eps |B_i|
\end{equation} 
for every $1 \leq i \leq k$. We will also ensure that whenever a vertex $v$ is moved to a cluster $A^*_i$, $v$ has at least $(d - \varepsilon)|B_i|$ neighbours in $B_i$, and vice versa.
We will terminate the procedure when $|A^*_i| = |A_i| + a_i$ and $|B^*_i| = |B_i| + b_i$, and then set $A'_i := A^*_i$ and $B'_i := B^*_i$ for each $1\leq i \leq k$.

Each iteration proceeds as follows: Let $1 \leq i \leq k$ be such that $|A^*_i| < |A_i| + a_i$ and let $j \neq i$ be such that $|A^*_j| > |A_j| + a_j$. 
(Such $i$ and $j$ exist by (\ref{iter1}).)
Suppose that $i < j$. 
Note that (i) implies that $(B_{j-1},A_j)_G$ is an $(\eps, d)$-regular pair. So
by (\ref{DifferenceBound}) and Fact \ref{RegularDegree} there is a vertex $v$ in $A^*_j$ which has at least $(d- \varepsilon)|B_{j-1}|$ neighbours in $B_{j-1}$. We move $v$ from $A^*_{j}$ to $A^*_{j-1}$. Similarly we move one vertex (which need not be $v$) from $A^*_{j-1}$ to $A^*_{j-2}$, and so on until we move one vertex from $A^* _{i+1}$ to  $A^*_i$. On the other hand, if $j < i$ we perform the same procedure moving vertices in the \emph{same} direction as before. That is, we move a vertex from $A^* _j$ to $A^* _{j-1}$ and
so on until we move a vertex $A^* _2 $ to $A^* _1$. Then we move a vertex $A^* _1$ to $A^* _k$ and continue until we move a 
vertex from $A^* _{i+1}$ to $A^* _i$.

Since in each step of the process we only move vertices between the $A^* _i$, certainly (\ref{iter1}) holds throughout.
Now when the procedure terminates we have $|A^*_i| = |A_i| + a_i$ for all $1 \leq i \leq k$. It remains to show that (\ref{DifferenceBound}) holds. Note that in each step of the iteration we add at most one vertex to each $A^*_i$ and remove at most one vertex from each $A^*_i$.
Further, in total we need to perform the iterative procedure at most 
$$\sum ^k  _{i=1} |a_i|+ \sum ^k _{i=1} a_i \leq (k+1)\xi n \leq 2 k \xi n$$ times. (The $\sum ^k _{i=1} a_i$ here comes from the
fact that, at the start, we moved $\sum ^k _{i=1} a_i$ vertices from $B_{i_2}$ to $A_{j_2}$.)
Thus, at the end of the procedure $|A^*_{j_2} \Delta A_{j_2}| \leq 5k\xi n$ and $|A^*_i \Delta A_i| \leq 4k\xi n$ for all $i \neq j_2$. We now set $A'_i := A^*_i$ for each $1 \leq i \leq k$.

We apply an identical iterative procedure to the $B^* _i$.
However, we now move vertices in the \emph{opposite} direction to before 
(so vertices are moved from $B^*_j$ to $B^*_{j+1}$, etc.).
Therefore we obtain sets $B'_i$ such that $|B'_i| = |B_i| + b_i$ and $|B'_i \Delta B_i| \leq 5k\xi n$ for all $1 \leq i \leq k$.

Given any $VW \in E(R)$, $(V,W)_G$ is an $(\eps , d)$-regular pair. Since by (\ref{DifferenceBound}),  $|V '\Delta V|$, $|W' \Delta 
W|\leq \eps |V|, \eps  |W|$, Proposition~\ref{PreserveReg} implies that $(V',W')_G$ is an $(\eps', d')$-regular pair. So indeed, $R$ is an $(\eps',d')$-reduced graph of $G$ on $A'_1, B'_1, A'_2, B'_2, \ldots, A'_k, B'_k$.
It remains to show that the pair $(A'_i, B'_i)_G$ is $(\varepsilon', d')$-super-regular for every $1 \leq i \leq k$. By (iv)
and (\ref{DifferenceBound}), every vertex $v \in A_i$ has at least $(d - \varepsilon)|B_i|\geq d'|B'_i|$ neighbours in $B'_i$. Further, during our iterative procedure we ensured that every vertex $v \in A'_i \backslash A_i$ has at least $(d - \varepsilon)|B_i|$ neighbours in $B_i$. Hence (\ref{DifferenceBound}) implies that every $v \in A'_i$ has at least
$$(d - \varepsilon)|B_i| - \eps |B_i|\geq  d'|B'_i|$$
neighbours in $B'_i$. Similarly each $w \in B'_i$ has at least $d'|A'_i|$ neighbours in $A'_i$. So $(A'_i, B'_i)_G$ is an $(\varepsilon', d')$-super-regular pair for all $1 \leq i \leq k$, as desired. \endproof
\COMMENT{Would a picture here help explain the movements made?}

\section{The Lemma for $G$} \label{SecLemmaForG}
\begin{lemma}[Lemma for $G$] \label{LemmaForG} Let $n_0 \in \mathbb{N}$ and let $\lambda, \xi, \varepsilon, d, \nu, \tau, \eta$ be positive constants such that 
$$0 < 1/n_0 \ll \lambda \ll \xi \ll \varepsilon \ll d \ll \nu \leq \tau \ll \eta \ll 1.$$
Suppose $G$ is a graph on $n \geq n_0$ vertices with $\delta(G) \geq \eta n$ which is a robust $(\nu, \tau)$-expander. Then there exists an integer $k$ such that $\xi \ll 1/k \ll \varepsilon$, 
integers $1\leq i_1 \not =j_1, i_2 \not = j_2 \leq k$
and a partition $(n_i)^{2k}_{i =1}$ of $n$ with $n_i > n/3k$ for all $1 \leq i \leq 2k$ and $|n_{2i-1} - n_{2i}| \leq \lambda n$ for all $1 \leq i \leq k$ such that the following holds: For every partition $(n'_i)^{2k}_{i =1}$ of $n$ satisfying $n'_i \leq n_i + \xi n$ for all $1 \leq i \leq 2k$, there exists a partition $A'_1, B'_1, A'_2, B'_2, \ldots, A'_k, B'_k$ of $V(G)$ and a spanning subgraph $G'$ of $G$ such that the following properties are satisfied.
\begin{itemize}
\item[\rm ($\alpha_1$)] $|A'_i| = n'_{2i-1}$ and $|B'_i| = n'_{2i}$ for all $1 \leq i \leq k$;
\item[\rm ($\alpha_2$)] $(A'_i, B'_i)_{G'}$ is $(\varepsilon, d)$-super-regular for all $1 \leq i \leq k$;
\item[\rm ($\alpha_3$)] $(B'_i, A'_{i+1})_{G'}$ is $(\varepsilon, d)$-regular for all $1 \leq i \leq k$ (where $A'_{k+1}: = A'_1$);
\item[\rm ($\alpha_4$)] $(A'_{i_1}, A'_{j_1})_{G'}$ is $(\varepsilon, d)$-regular;
\item[\rm ($\alpha_5$)] $(B'_{i_2}, B'_{j_2})_{G'}$ is $(\varepsilon, d)$-regular.
\end{itemize}
\end{lemma}

\proof Choose additional constants $\varepsilon'$ and $d'$ such that
$$\xi \ll \varepsilon' \ll \varepsilon \ll d \ll d' \ll \nu.$$
Apply the Regularity Lemma (Lemma \ref{Regularity}) with parameters $\varepsilon'$, $d'$ and $k_0 := 1/\varepsilon'$ to obtain clusters $V_1, \ldots, V_{k'}$ of size $m$  (where $(1-\varepsilon')n/k' \leq m \leq n/k'$), an exceptional set $V_0$, a pure graph $G' \subseteq G$ and the reduced graph $R$ of $G$ with parameters $\eps ' , d'$ and $k_0$. Since $\xi \ll \varepsilon '$ we may assume that
$$\xi \ll 1/k' \leq \eps '.$$

If $k'$ is odd then we delete $V_{k'}$ from $R$ and add all of the vertices of $V_{k'}$ to $V_0$. So $|V_0| \leq \varepsilon' n + m \leq 2 \varepsilon' n$. We now refer to this modified reduced graph as $R$ and redefine $k' = |R|$. By Lemma~\ref{PreserveRE}, $R$ originally had minimum degree at least $\eta k'/2$ and was a robust $(\nu/2, 2\tau)$-expander. So $R$ still has minimum degree at least $\eta k'/3$ and by Lemma \ref{PreserveRobustExp}, $R$ is still a robust $(\nu/3, 3\tau)$-expander. 

Set $k := k'/2$.
Since $1/k' \ll \nu \leq \tau \ll \eta < 1$, Theorem \ref{expanderthm} implies that $R$ contains a Hamilton cycle $C = A_1 B_1 \ldots A_k B_k A_1$.  Since $|C| = 2k$ is even, $C$ contains a perfect matching $M = \{A_1 B_1, \ldots, A_k B_k\}$.
Notice that $R$ contains an edge $A _{i_1} A _{j_1}$ for some $1\leq i_1 \neq j_1 \leq k$ and an edge $B _{i_2} B _{j_2}$ for some
 $1\leq i_2 \neq j_2\leq k$. Indeed, let $A := \{A_i\}_{i = 1}^k$ and note that since $R$ is a robust $(\nu/3, 3\tau)$-expander we have $|RN_{\nu, R}(A)| \geq k + \nu k'$. This implies that $A \cap RN_{\nu, R}(A) \neq \emptyset$ and hence that $R$ contains some edge $A_{i_1} A_{j_1}$. Similarly $R$ contains an edge $B_{i_2} B_{j_2}$.

Fact \ref{RegularDegree} implies that we can replace each cluster in $ V(R)$ with a subcluster of size $m':=(1-\eps ')m$
such that for every edge $A_jB_j \in M$ the chosen subclusters of $A_j$ and $B_j$ form a $(2\eps' , d'/2)$-super-regular pair in
$G'$.  We add all of the vertices not in these subclusters to $V_0$, and from now on we refer to the subclusters as the clusters of $R$.
So $(V,W)_{G'}$ is still a $(2\eps ', d'/2)$-regular pair for all $VW \in E(R)$.
Note that $|V_0| \leq 2\varepsilon'n + \varepsilon' n = 3 \varepsilon' n$.

Our next task is to incorporate the vertices of $V_0$ into the clusters $V_1, \ldots, V_{k'}$ such that the pairs $(A_j, B_j)_{G'}$
remain super-regular and such that the pairs $(V_i, V_j)_{G'}$ remain regular for all $V_iV_j \in E(R)$ (with somewhat weaker constants in each case). Let $V_0 = \{x_1, \ldots, x_t\}$ where $t \leq 3 \varepsilon' n$. We will assign the vertices of $V_0$ in such a way that:
\begin{itemize}
\item[\rm (a)] At most $8 \varepsilon' m'/\eta$ vertices are assigned to each cluster $V \in V(R)$;
\item[\rm (b)] Whenever a vertex $x_i \in V_0$ is assigned to a cluster $A_j$, $x_i$ has at least $\eta m'/4$ neighbours in $B_j$. Similarly any vertex from $ V_0$ assigned to $B_j$ has at least $\eta m'/4$ neighbours in $A_j$.
\end{itemize} 
Suppose we have assigned $x_1, \ldots, x_{i-1}$ to clusters in $V(R)$ such that (a) and (b) are satisfied. Call a cluster $V \in V(R)$ \emph{full} if it has already been assigned $8 \varepsilon' m'/\eta$ vertices of $V_0$. Let $F$ be the set of full clusters. Since $|V_0| \leq 3 \varepsilon' n$ we have $|F| \leq (3 \varepsilon' \eta n)/(8 \varepsilon' m') \leq \eta k$. Thus, as $\delta(G) \geq \eta n$,
$$\left|N_G(x_i) \backslash \left(V_0 \cup \bigcup_{V \in F} V\right)\right| \geq \eta n - 3 \varepsilon' n - (\eta k)m' \geq \eta n/3.$$
Hence, by the pigeonhole principle there exists some $V \in V(R) \backslash F$ such that $|N_G(x_i) \cap V| \geq \eta n/3k' \geq \eta m'/4$. Now if $V = A_j$ for some $1 \leq j \leq k$ then we add $x_i$ to $B_j$; otherwise, $V = B_j$ for some $1 \leq j \leq k$ and we add $x_i$ to $A_j$. Repeating this process for each $x_i$ we indeed assign all of the vertices of $V_0$ to the clusters of $R$ in such a way that (a) and (b) are satisfied. We now incorporate all of the assigned vertices into their respective clusters. 
Further, we add all those edges from $G$ with endpoints in $V_0$ to $G'$.
Note that 
\begin{itemize}
\item[\rm (c)] $m' \leq |V| \leq m' + 8\varepsilon' m'/\eta \leq (1 + \sqrt{\varepsilon'})m'$ for all $V \in V(R)$,
\item[\rm (d)] $(V, W)_{G'}$ is a $((\varepsilon')^{1/3}, d'/4)$-regular pair for every edge $VW \in E(R)$ and
\item[\rm (e)] $(V, W)_{G'}$ is a $((\varepsilon')^{1/3}, d'/4)$-super-regular pair for every edge $VW \in E(M)$.
\end{itemize}
(Conditions (d) and (e) follow by Proposition~\ref{PreserveReg}.)

Next we will perform an algorithm which redistributes vertices among the clusters in $R$ in such a way that $||A_i| - |B_i|| \leq \lambda n$ for each $1 \leq i \leq k$. We define $ \{A^*_i, B^*_i\}_{i = 1}^k$, $R^*$ and $M^*$ as follows: Initially 
we set $A^*_i := A_i$ and $B^*_i := B_i$ for all $1 \leq i \leq k$, $R^* := R$ and $M^* := M$. At each step we will redefine each $A^*_i$ and $B^*_i$, $R^*$ and $M^*$ and reassign vertices so that the quantity
$$\Sigma^* = \sum_{1 \leq i \leq k, ||A^*_i| - |B^*_i|| > \lambda n} ||A^*_i| - |B^*_i||$$
decreases by at least $\lambda n$.  The algorithm will terminate when $\Sigma^* = 0$, i.e. when $||A^*_i| - |B^*_i|| \leq \lambda n$ for all $1 \leq i \leq k$.  Initially $\Sigma^* \leq 8 \varepsilon' m' k/\eta \leq 4 \varepsilon' n/\eta$ by (c), and hence we need at most $4 \varepsilon' /\eta \lambda$ steps to complete the process. $R^*$ will always be an \emph{induced}
 subgraph of $R$ and at each step we set $M^* $ to be the submatching of $M$ induced by $V(R^*)$.
 (Note that $V(R^*)$ is a subset of $V(R)=\{A_i, B_i\}_{i = 1}^k$ throughout the algorithm.)

We will ensure that the inequality
\begin{equation} \label{LFGIneq1}
|R^*| \geq \left(1 - \nu/12\right)k'
\end{equation} 
holds throughout, and that $M^*$ is a perfect matching in $R^*$. Further we will ensure that
\begin{equation} \label{LFGIneq2}
|A^*_i \backslash A_i| \leq (\varepsilon')^{1/3} m' \text{ and } |B^*_i \backslash B_i| \leq (\varepsilon')^{1/3} m';
\end{equation}
\begin{equation} \label{LFGIneq3}
|A_i \backslash A^*_i| \leq (\varepsilon')^{1/3} m' \text{ and } |B_i \backslash B^*_i| \leq (\varepsilon')^{1/3} m'
\end{equation}
for all $1 \leq i \leq k$.

Each step proceeds as follows: Call a vertex $v$ \emph{well-connected} to a cluster $V \in V(R)$ if $v$ has at least $d' m'/8$ neighbours in $V$. Recall that if $VW\in E(R)$ then $(V, W)_{G'}$ is a $((\varepsilon')^{1/3}, d'/4)$-regular pair and so $V$ contains at least $m'/2$ vertices $v$ which are well-connected to $W$. In what follows we will ensure that every vertex we redistribute to a cluster $A^*_i$ is well-connected to $B_i$ and vice versa. 
Since (\ref{LFGIneq3}) holds throughout the process, 
given any $VW \in E(R^*)$, 
$V^*$ will always contain at least $m'/2-(\varepsilon')^{1/3} m'\geq m'/3 \gg \lambda n$
vertices that are well-connected to $W$ (where $V^*:=A^*_i$ if $V=A_i$ for some $i$ and $V^*:=B^*_i$ if $V=B_i$ for some $i$).  Thus, at any point during the algorithm we may choose a set of $\lambda n$ well-connected vertices from any of the $A^*_i$ and $B^*_i$.
(When it is clear from the context, we will not explicitly specify which cluster a vertex $v$ is well-connected to.)

Let $S$ be the set of clusters $V^* \in V(R^*)$ such that either $V^* = A_i$ where $|A^*_i| > |B^*_i| + \lambda n$ or $V^* = B_i$ where $|B^*_i| >|A^*_i| + \lambda n$. If $S$ is empty then the algorithm terminates. (We shall see later that in this case we must
have that $\Sigma ^*=0$.) Otherwise, choose $V^* \in S$ arbitrarily. Since $R$ is a robust $(\nu/3, 3\tau)$-expander and $\delta(R) \geq \eta k'/3$, (\ref{LFGIneq1}) implies that $\delta(R^*) \geq \eta |R^*|/4$ and Lemma \ref{PreserveRobustExp} implies that $R^*$ is a robust $(\nu/4, 4 \tau)$-expander. Hence Lemma \ref{FindWalk} implies that $R^*$ contains a shifted $M^*$-walk $P'$ of length at most $12/\nu$ which starts and finishes at $V^*$.
By Lemma \ref{FindSimpleWalk}, $P'$ contains a simple shifted $M^*$-walk $P''$ which also starts and finishes at $V^*$. Now apply Lemma \ref{FindPureWalk} to $P''$ to obtain a simple shifted $M^*$-walk $P$ of length at most $12/\nu$, such that the endpoints of $P$ both lie in $S$ and no other vertices of $P$ lie in $S$. We call $P$ the \emph{active walk} of this step of the algorithm.

Let $P = U_1 W_2 U_2 \ldots W_{\ell-1} U_{\ell-1} W_{\ell}$ such that $W_i U_i \in E(M^*)$ for each $2 \leq i \leq \ell-1$. Let $W_1$ and $U_{\ell}$ denote the clusters such that $W_1 U_1, W_\ell U_\ell \in E(M^*)$. 
Given any $1 \leq i \leq \ell$, if $U_i=A_j$ for some $1 \leq j \leq k$, set $U^* _i:=A^*_j$; otherwise $U_i=B_j$ for some $1 \leq j \leq k$, so set $U^* _i:=B^*_j$. Define $W^* _i$ analogously for each $1 \leq i \leq \ell$.
Move $\lambda n/2$ well-connected vertices from $U^*_1$ into $U^*_2$, $\lambda n/2$ well-connected vertices from $U^* _2$ into $U^* _3$, and so on until we have moved $\lambda n/2$
well-connected vertices from $U^* _{\ell -1}$ to $U^*_\ell$. Then move $\lambda n/2$ well-connected vertices from $W^*_2$ into $W^*_1$, $\lambda n/2$ well-connected vertices from $W^*_3$ into $W^*_2$, and so on until we have moved $\lambda n/2$
well-connected vertices from $W^*_{\ell}$ to $W^*_{\ell-1}$.
Note that since $P$ is simple, each cluster loses at most $\lambda n$ vertices
and gains at most $\lambda n$ vertices. Further for each $1 < i < \ell$ the quantity $||W^*_i| - |U^*_i||$ remains unchanged
(in fact, $|W^*_i|$ and $|U^*_i|$ remain unchanged). 
For $i = 1, \ell$, $||W^*_i| - |U^*_i||$ decreases by precisely $\lambda n$ (or $2 \lambda n$ if $U_1 = W_\ell$).

In order to ensure that (\ref{LFGIneq2}) holds, we remove from $R^*$ every pair $\{A_i, B_i\}$ of clusters such that $|A^*_i \backslash A_i| \geq (\varepsilon')^{1/3} m' - \lambda n$ or $|B^*_i \backslash B_i| \geq (\varepsilon')^{1/3} m' - \lambda n$. Since each cluster gains at most $\lambda n$ vertices in each step, any clusters which are not removed will still satisfy (\ref{LFGIneq2}) at the end of the next step. Similarly,  to ensure that (\ref{LFGIneq3}) holds, we remove from $R^*$ every pair $\{A_i, B_i\}$ of clusters such that $|A_i \backslash A^*_i| \geq (\varepsilon')^{1/3} m' - \lambda n$ or $|B_i \backslash B^*_i| \geq (\varepsilon')^{1/3} m' - \lambda n$.  \medskip

\begin{claim}\label{claim1}
 For each pair $\{A_i, B_i\}$ of clusters which is removed from $R^*$ we have that $||A^*_i| - |B^*_i|| \leq \lambda n$.
\end{claim}
\noindent To prove the claim, suppose for a contradiction that some pair $\{A_i, B_i\}$ of clusters is removed from $R^*$ and that $||A^*_i| - |B^*_i|| > \lambda n$. In order for $\{A_i, B_i\}$ to be removed we must have that $|A^*_i \backslash A_i| \geq (\varepsilon')^{1/3} m' - \lambda n$, $|B^*_i \backslash B_i| \geq (\varepsilon')^{1/3} m' - \lambda n$, 
$|A_i \backslash A^*_i| \geq (\varepsilon')^{1/3} m' - \lambda n$ or $|B_i \backslash B^*_i| \geq (\varepsilon')^{1/3} m' - \lambda n$. Without loss of generality assume that $|A^*_i \backslash A_i| \geq (\varepsilon')^{1/3} m' - \lambda n$. Since in each step we add at most $\lambda n$ vertices to $A^*_i$, there must have been at least $(\varepsilon')^{1/3} m'/2 \lambda n$ steps in the algorithm so far, such that $A_i$ is contained in the active walk $P$ of each step. By the definition of $P$, either $A_i$ or $B_i$ must be an endpoint of $P$. So $||A^*_i| - |B^*_i||$ is reduced by at least $\lambda n$ during each such step, and hence by at least $(\varepsilon')^{1/3} m'/2$ during the algorithm so far. But this is a contradiction since initially $||A^*_i| - |B^*_i|| \leq \sqrt{\varepsilon'} m'\leq (\eps ')^{1/3} m'/2$. \medskip

It remains to show that (\ref{LFGIneq1}) holds throughout the process. Suppose for a contradiction that at some point more than $\nu k'/12$ clusters have been removed from $R^*$. Then at least 
\COMMENT{I think this is actually correct. Every time we move $(\varepsilon')^{1/3} m' - \lambda n$ vertices from $A^*_i$ to $A^*_j$ we remove four clusters ($A^*_i$, $A^*_j$, $B^*_i$, $B^*_j$), so we need to divide by four. Still should we put the $1/4$ outside the bracket?}
\begin{equation} \label{LFGIneq4} 
((\varepsilon')^{1/3} m' - \lambda n)\left(\frac{\nu k'}{48}\right) > \sqrt{\varepsilon'} n 
\end{equation}
vertices of $G$ must have been redistributed during the process so far. But at most $12 \lambda n/\nu$ vertices were redistributed during each step, and at most $4 \varepsilon' /\eta \lambda$ steps were performed during the process. Hence the number of redistributed vertices is at most 
$$\frac{12 \lambda n}{\nu} \cdot \frac{4 \varepsilon'}{\eta \lambda} < \sqrt{\varepsilon'} n,$$
which contradicts (\ref{LFGIneq4}). This proves that (\ref{LFGIneq1}) holds throughout. 

By construction, when the algorithm terminates we have that $V(R^*)$ does not contain any $A_i$ and $B_i$ such that
$||A^*_i|-|B^* _i||> \lambda n$. Further,  by Claim~\ref{claim1}, for those clusters $A_i,B_i \notin V(R^*)$
we have that $||A^*_i|-|B^* _i|| \leq \lambda n$. Hence, we indeed obtain clusters $ \{A^*_i, B^*_i\}_{i = 1}^k$
such that $\Sigma ^* =0$ and (\ref{LFGIneq2}) and (\ref{LFGIneq3}) hold.

We now set $n_{2i - 1}: = |A^*_i|$ and $n_{2i} := |B^*_i|$ for each $1 \leq i \leq k$. Notice that $n_j \geq (1 - (\varepsilon')^{1/3})m' > n/3k$ for each $1 \leq j \leq 2k$. 
We now relabel the clusters of $R$ in the natural way so that $V(R) = \{A^*_i, B^*_i\}_{i = 1}^k$. 
Note that by (\ref{LFGIneq2}) and (\ref{LFGIneq3}) we have
$$|A_i \Delta A^*_i| \leq 2(\varepsilon')^{1/3} m' \text{ and } |B_i \Delta B^*_i| \leq 2(\varepsilon')^{1/3} m'$$
for each $1 \leq i \leq k$. Hence by Proposition \ref{PreserveReg}  the pair $(V, W)_{G'}$ is $((\varepsilon')^{1/10}, d'/10)$-regular for every edge $VW\in E(R)$. Further, the pair $(A^*_i, B^*_i)_{G'}$ is $((\varepsilon')^{1/10}, d'/10)$-super-regular for every $1 \leq i \leq k$. Indeed, we ensured that every vertex $v$ which was redistributed to $A^*_i$ had at least $d'm'/8$ neighbours in $B_i$. Since $|B^*_i \Delta B_i| \le 2 (\varepsilon')^{1/3} m'$,
$v$ has at least
$d'|B^*_i|/10$ neighbours in $B^*_i$. 
Similarly every vertex $w \in B^*_i$ has at least $d'|A^*_i|/10$ neighbours in $A^*_i$.

Now suppose we are given a partition $(n'_i)_{i = 1}^{2k}$ of $n$ such that $n'_i \leq n_i + \xi n$ for each $1 \leq i \leq 2k$. 
Set $a_i := n'_{2i-1} - n_{2i-1}$ and $b_i := n'_{2i} - n_{2i}$ for all $1\leq i \leq k$. 
Notice that $|a_i|, |b_i| \leq 2k\xi n$ for each $1 \leq i \leq k$ and $|\sum_{i = 1}^k a_i| = |\sum_{i = 1}^k b_i| \leq 2k\xi n$. 
Recall that $R$ contains the edges $A^* _{i_1}A^*_{j_1}$ and $B^*_{i_2}B^*_{j_2}$.
  Thus, we can apply the Mobility lemma (Lemma \ref{MobLemma}) with parameters $k, 2k\xi, (\varepsilon')^{1/10}, \varepsilon, d$ and $d'/10$ to obtain a partition $A'_1, B'_1, \dots ,A'_k,B'_k$ of $V(G)$ which satisfies conditions ($\alpha_1$)--($\alpha_5$).\COMMENT{Include $*$ since we relabeled $V(R)=V(R)^*$.}
 \endproof

\section{The Lemma for $H$} \label{SecLemmaForH}

\begin{lemma}[Lemma for $H$] \label{LemmaForH} For any $\Delta, k \in \mathbb{N}$ and $\xi > 0$, there exist $\beta > 0$ and $n_0 \in \mathbb{N}$ such that the following holds: Let $H$ be a bipartite graph on $n\geq n_0$ vertices with bandwidth at most $\beta n$ and such that $\Delta(H) \leq \Delta$. Let $n_1, n_2, \ldots, n_{2k}$ be an integer partition of $n$ such that $n_i >n/(3k)$ for all $1\leq i
\leq 2k$ and 
 $|n_{2i-1} - n_{2i}| \ll \xi n$ for all $1 \leq i \leq k$. Suppose $C$ is the cycle $12 \ldots (2k)1$ on $[2k]$, and let $c = \{2i_1, 2i_2\}$  be a chord of $C$ (for some distinct $1\leq i_1 , i_2 \leq k$). Then there exists a set $S \subseteq V(H)$ and a graph homomorphism $f:H \rightarrow C \cup \{c\}$, such that
\begin{itemize}
\item[($\beta _1$)] $|S| \leq \xi n$;
\item[($\beta _2$)] $|f^{-1}(i)| \leq n_i + \xi n$ for all $1 \leq i \leq 2k$;
\item[($\beta _3$)] Every edge which is not in $H[S]$ is mapped to an edge $\{2i-1,  2i\}$, for some $1 \leq i \leq k$.
\end{itemize}
\end{lemma}

\proof  Choose $\beta > 0$ and integers $n_0$, $m_1$, $m_2$ and $k_1$ such that
$$1/n_0 \ll \beta \ll 1/m_1 \ll 1/m_2 \ll 1/k_1 \ll 1/\Delta, 1/k, \xi.$$
Further, we may assume that $m_2$ divides $m_1$.
We begin by defining a new cycle $C'$ with chord $c'$ which will act as an intermediate stage between $C$ and $H$, i.e., we will construct homomorphisms $f_1: C' \cup \{c'\} \rightarrow C \cup \{c\}$ and $f_2: H \rightarrow C' \cup \{c'\}$ such that $f = f_1 \circ f_2$ is our desired homomorphism. The homomorphism $f_2$ will be constructed to map roughly the same number of vertices to 
each vertex in $C'$. Notice however, that our desired homorphism $f$ may not map vertices in an `equal' way (since, in general,  the
$n_i$ may be  far from equal). Thus, the role of $f_1$ is to ensure $f$ maps the `correct' proportion of vertices
to each vertex in $C$.

Let $C'$ be the cycle $12 \ldots (2k')1$ on $[2k']$, where $k': = \sum_{i \in [k]} \lceil (n_{2i-1} + n_{2i}) k_1/n \rceil$. Note that $ k_1 \leq k' \leq  k_1 + k$. We define $f_1$ as follows: For each $1 \leq j \leq k'$, let $g(j)\in \mathbb N$ be such that $\sum_{i = 1}^{g(j)-1} \lceil (n_{2i-1} + n_{2i}) k_1/n \rceil < j \leq \sum_{i = 1}^{g(j)} \lceil (n_{2i-1} + n_{2i}) k_1/n \rceil$. Then set $f_1(2j-1) = 2g(j) - 1$ and $f_1(2j) = 2g(j)$ for each $1 \leq j \leq k'$. 

Recall that $c=\{2i_1, 2i_2\}$ is a chord of $C$. Suppose that $i'_1, i'_2$ are such that $f_1 (2i'_1)=2i_1$ and $f_1 (2i'_2)=2i_2$.
Notice that as $i_1 \not = i_2$, we have that $i'_1 \not = i'_2$. Thus, set $c':=\{2i'_1, 2i'_2\}$ to be the chord of $C'$.

By construction $f_1 (c')=c$. Given any edge $c_1=\{2j-1,2j\}$ on $C'$, we have that $f_1 (c_1)=\{2g(j)-1, 2g(j)\}$. Further,
consider any edge $c_2=\{2j,2j+1\}=\{2j, 2(j+1)-1\}$ on $C'$. Then $f_1(2j)=2g(j)$ and $f_1 (2(j+1)-1)=2g(j+1)-1$. By definition of $f_1$,
either $g(j+1)=g(j)$ or $g(j+1)=g(j)+1$. But both $\{2g(j),2g(j)-1 \}$ and $\{2g(j),
2g(j)+1\}$ are edges of $C$. So in either case $f_1$ maps
$c_2$ to an edge of $C$. Therefore, indeed $f_1$ is a graph homomorphism.

Roughly speaking, we will construct $f_2$ as follows: Initially we split $H$ up into \emph{small segments} $A_1, B_1, \ldots, A_{m_1}, B_{m_1}$ in such a way that almost all of the edges of $H$ lie in the pairs $(A_i, B_i)_{i = 1}^{m_1}$ and the remainder lie in the pairs $(B_i, A_{i+1})_{i = 1}^{m_1}$. Our ideal strategy would be to map all of the vertices of $A_1$ onto vertex $1$ of $C'$, the vertices of $B_1$ onto vertex $2$, the vertices of $A_2$ onto vertex $3$, etc. This ensures that $f_2$ is a homomorphism and that almost all of the edges of $H$ are mapped onto an edge of the form $\{2i-1, 2i\}$ for some $i$. However the number of vertices mapped onto each vertex of $C'$ may vary widely. To solve this problem we introduce `drunken' segments in which the assignment of the vertices is random, and use a probabilistic argument to show that with positive probability each vertex of $C'$ receives approximately the same number of vertices of $H$. We also use  the chord $c'$  to `turn around' at some point during the process, in order to eliminate the possible inequality between the number of vertices of $H$ assigned to odd and even vertices of $C'$.
\medskip

\noindent
{\bf Chopping $H$ up into segments.}
Since $H$ has bandwidth at most $\beta n$, there exists an ordering $x_1, x_2, \ldots, x_n$ of $V(H)$ such that for every edge $x_i x_j$ of $H$, $|i - j| \leq \beta n$. Let $(A, B)$ be a bipartition of $V(H)$. We define $\{A_i, B_i\}_{i = 1}^{m_1}$ as follows: for each vertex $x_s \in A$ there exists $1 \leq i\leq m_1$ such that $(i-1)n/m_1 - \beta n < s \leq in/m_1 - \beta n$ (unless $s > n - \beta n$). We assign $x_s$ to $A_i$ (or to $A_{m_1}$ if $s > n - \beta n$). Similarly for each vertex $x_t \in B$ there exists $1\leq j \leq m_1$ such that $(j-1)n/m_1 < t \leq jn/m_1$, and we assign $x_t$ to $B_j$. Let $S$ be the set of vertices $x_s$ such that $in/m_1 - 2\beta n < s \leq in/m_1 + \beta n$ for some $1 \leq i \leq m_1$. Note that the following properties hold:
\begin{itemize}
\item[\rm (a)] $n/m_1 - \beta n \leq |A_i| + |B_i| \leq n/m_1 + \beta n$ for each $1 \leq i \leq m_1$;
\item[\rm (b)] $|S| \leq 3m_1 \beta n \ll \xi n$;
\item[\rm (c)] Every edge of $H$ which is not in $H[S]$ lies in one of the pairs $(A_i, B_i)$ for some $1 \leq i \leq m_1$;
\item[\rm (d)] Every edge of $H[S]$ lies  in one of the pairs $(A_i, B_i)$ or one of the pairs $(B_i, A_{i+1})$ for some $1 \leq i \leq m_1$ (where $A_{m_1+1}:=A_1$).
\end{itemize}
Properties (c) and (d) follow from the fact that $H$ has bandwidth at most $\beta n$
and that $n/m_1 \gg \beta n$. We now modify the small segments so that properties (a)--(d) are still satisfied and so that every small segment has size at least $n/(4 \Delta m_1)$. Suppose a small segment $A_i$ has size smaller than $n/(4 \Delta m_1)$. Note that $|N_H(A_i)| \leq n/(4 m_1)$ and so 
$$|B_i \backslash (S \cup N_H(A_i))| \stackrel{(a)}{\geq} (n/m_1 - \beta n - n/(4 \Delta m_1)) - 3 \beta n - n/(4m_1) \geq n/(4m_1).$$ 
But (c) implies that any vertex in $B_i \backslash (S \cup N_H(A_i))$ must be isolated in $H$ and so may be reassigned to $A_i$ without affecting properties (a)--(d). Hence we may reassign sufficiently many vertices so that $|A_i|, |B_i| \geq n/(4 \Delta m_1)$. For any segment $B_i$ which has size smaller than $n/(4 \Delta m_1)$ we proceed in an identical way. From now on we denote by $A$ the union of small segments $\bigcup_{i = 1}^{m_1} A_i$ and by $B$ the union $\bigcup_{i = 1}^{m_1} B_i$.

We now group the small segments together to form \emph{large segments} $\{L_j\}_{j = 1}^{m_2}$, which are defined as
$$L_j := \bigcup_{\frac{(j-1)m_1}{m_2} < t \leq \frac{jm_1}{m_2}} \left(A_t \cup B_t\right).$$
Note that since $\beta \ll 1/m_1$, (a) implies that
\begin{equation} \label{L4H1}
\frac{n}{m_2} - \sqrt{\beta} n \leq |L_j| \leq \frac{n}{m_2} + \sqrt{\beta} n
\end{equation}
for each $1 \leq j \leq m_2$. In order to eliminate any inequality between the number of vertices of $H$ assigned to odd and even vertices of $C'$ we need to partition $\{L_j\}_{j = 1}^{m_2}$ into two parts. We will assign the vertices in each part separately. For each $1 \leq j \leq m_2$, set $s_j := |L_j \cap A| - |L_j \cap B|$. Note that $|s_j| \leq n/m_2 + \sqrt{\beta} n - n/(4\Delta m_2) \leq n/m_2$ for each $1 \leq j \leq m_2$, and that $\sum_{j = 1}^{m_2} s_j = |A| - |B|$. Suppose without loss of generality that $|A| - |B| \geq 0$. Then since $\beta, 1/m_2 \ll \xi$ there exists an integer $m_3$ so that  $\xi m_2/20 \leq m_3 \leq (1 - \xi/20)m_2$ and 
\begin{equation} \label{L4H2}
\frac{(|A| - |B|)}{2} - \frac{\xi n}{20} \leq \sum_{i = 1}^{m_3} s_j \leq \frac{(|A| - |B|)}{2} + \frac{\xi n}{20}.
\end{equation}
We will embed separately the large segments $\{L_j\}_{j = 1}^{m_3}$ and the segments $\{L_j\}_{j = m_3 + 1}^{m_2}$.

For each $1 \leq j \leq m_3$, let the \emph{drunken segment} $D_j$ be the union of the last $k_2 := \xi m_1/(6 k' m_2)$ pairs of small segments in $L_j$ and let the \emph{sober segment} $S_j$ be the union of the rest of the small segments in $L_j$. 
\medskip

\noindent
{\bf Defining our algorithms.}
We now define three different algorithms for assigning the vertices of a segment to vertices of $C'$, given an \emph{initial vertex} $2i'_0-1 \in C'$. In each case we work mod $2k'$ when dealing with vertices of $C'$. The \emph{sober algorithm} is a deterministic process which proceeds as follows: 
Let $S_j$ be a sober segment whose first small segments are $A_{i_0}$ and $B_{i_0}$. For every $i$ such that $A_i$ and $B_i$ are small
segments of $S_j$, assign every vertex of $A_i$ to the vertex $2i'-1$ of $C'$ and every vertex of $B_i$ to the vertex $2i'$ of $C'$
 where $i' \equiv i'_0 + i - i_0 \mod k'$. So whenever $B_i$ is assigned to $2i'$, $A_{i+1}$ is assigned to $2i'+1$ ($\bmod$ $2k'$). We call the vertex $2i^*$ of $C'$ to which the vertices of the last small segment of $S_j$ are assigned the \emph{final vertex} of the algorithm and define this term in a similar way for the remaining two algorithms.

The \emph{drunken algorithm} is a randomised algorithm which proceeds as follows: Given a drunken segment $D_j$ whose first small segment is $A_{i_0}$, assign every vertex of $A_{i_0}$ to the vertex $2i'_0 - 1$ of $C'$ and every vertex of $B_{i_0}$ to the vertex $2i'_0$ of $C'$. Then for every pair $A_{i+1}, B_{i+1}$ of small segments in $D_j$, let $2i'$ be the vertex to which the vertices of $B_i$ were assigned and let
$$i'' = \begin{cases}
i' &\text{ with probability } \frac{1}{2}; \\
i' + 1 &\text{ with probability } \frac{1}{2}.
\end{cases}$$
(All random choices are made independently.) Assign every vertex of $A_{i+1}$ to $2i''-1$ and every vertex of $B_{i+1}$ to $2i''$. \medskip

\begin{claim}\label{claima} Suppose that the vertices of $D_j$ are assigned using the drunken algorithm with initial vertex $2i'_0 - 1$, and let $i'_1 \in [k']$ be arbitrary. Let the random variable $I$ be the final vertex of the drunken algorithm. Then
$$\mathbb{P}[I = 2i'_1 \mid i'_0] \leq \frac{1 + \xi/20}{k'}.$$ 
\end{claim}
\noindent To prove the claim, note that $I \sim 2(i'_0 + Bin(k_2, 1/2))$, that $k_2 \gg (k')^3/6$ and that $1/k' \ll \xi/20$. So Lemma \ref{EvenBinomial} with $\varepsilon := \xi/20$ implies that $\mathbb{P}[I - 2i'_0 = 2i'_1 - 2i'_0] \leq (1 + \xi/20)/k'$, and Claim \ref{claima} follows immediately. \medskip

The $2i'_1$-\emph{seeking algorithm} is a deterministic algorithm which proceeds as follows: Given a drunken segment $D_j$ whose first small segment is $A_{i_0}$, assign every vertex of $A_{i_0}$ to the vertex $2i'_0 - 1$ of $C'$ and every vertex of $B_{i_0}$ to the vertex $2i'_0$ of $C'$. Then for every pair $A_{i+1}, B_{i+1}$ of small segments in $D_j$, let $2i'$ be the vertex of $C'$ to which the vertices of $B_i$ were assigned and let
$$i'' = \begin{cases}
i' &\text{ if } i' = i_1; \\
i' + 1 &\text{ otherwise.}
\end{cases}$$
Assign every vertex of $A_{i+1}$ to $2i''-1$ and every vertex of $B_{i+1}$ to $2i''$. Note that the final vertex of this algorithm is always $2i'_1$, since $k' \leq \xi m_1/(6k'm_2)$.
\medskip

\noindent 
{\bf Applying the algorithms.}
We use these algorithms to assign small segments to vertices of $C'$ as follows: Choose $i'_0 \in [k']$ randomly and let $2i'_0 - 1$ be the initial vertex for $S_1$. For each $1 \leq j \leq m_3 - 1$, use the sober algorithm to assign the segments of $S_j$ and then use the drunken algorithm to assign the vertices of $D_j$, where in each case the initial vertex of each segment is the successor of the final vertex of the previous segment. 
(So for example, if the final vertex of the drunken algorithm, when applied to $D_j$, is $2i^*$, then the initial vertex of the
sober algorithm, when applied to $S_{j+1}$, is $2i^*+1$.)
Then assign the vertices of $S_{m_3}$ using the sober algorithm and assign the vertices of $D_{m_3}$ using the $2i'_1$-seeking algorithm. (Recall that $2i'_1$ was a vertex of the chord $c'$.) We explain how we assign the small segments from $\bigcup_{j = m_3 + 1}^{m_2} L_j$
later.

\begin{claim}\label{claimb}
For each $1 \leq i \leq k'$, let $X_i$ be the number of vertices of $\bigcup_{j = 1}^{m_3} S_j$ assigned to the vertex $2i-1$ of $C'$. Then
$$\mathbb{P}\left[X_i > \frac{1}{2k'}\left(\frac{m_3 n}{m_2} + \frac{|A| - |B|}{2}\right) + \frac{\xi n}{6k'}\right] \leq \frac{1}{3k'}.$$
\end{claim}
\noindent For each $1 \leq j \leq m_3$, let $Y_j = |S_j \cap A|$ and let $X_{i, j}$ be the number of vertices of $S_j$ which are assigned to $2i-1$. To prove the claim, we first use Claim \ref{claima} to bound $\mathbb{E}[X_{i, j} \mid X_{i, j-1}, \ldots, X_{i, 1}]$. Let $r_j$ be the initial vertex of $S_j$ for each $j$. Let $B$ be the event $X_{i, j-1} = x_{i, j-1}, \ldots, X_{i, 1} = x_{i, 1}$ for some $x_{i, 1}, \ldots x_{i, j-1}$. Now for $1 \leq i' \leq k'$ and any integer $x $ we have 
$\mathbb P [X_{i,j} =x \mid B \cap (r_{j-1} = 2i'-1)]= \mathbb P [X_{i,j} =x \mid  r_{j-1} = 2i'-1]$.
 Hence Lemma \ref{ConditionalExp} implies that $\mathbb{E}[X_{i, j} \mid X_{i, j-1}, \ldots, X_{i, 1}] \leq \max_{i' = 1}^{k'} \mathbb{E}[X_{i, j} \mid r_{j-1} = 2i'-1]$. But Claim \ref{claima} implies that
\begin{align*}
\mathbb{E}[X_{i, j} \mid r_{j-1} = 2i'-1] &= \sum_{i'' = 1}^{k'} \mathbb{E}[X_{i, j} \mid r_j = 2i''-1] \mathbb{P}[r_{j} = 2i''-1 \mid r_{j-1} = 2i'-1] \\
&\leq \frac{1 + \xi/20}{k'} \sum_{i'' = 1}^{k'} \mathbb{E}[X_{i, j} \mid r_j = 2i''-1] = \frac{(1 + \xi/20)Y_j}{k'}.
\end{align*}
Hence $\mathbb{E}[X_{i, j} \mid X_{i, j-1}, \ldots, X_{i, 1}] \leq (1 + \xi/20)Y_j/k'$. Set $X'_{i, j} := X_{i, j} m_2/n$. 
Since
$$|S_j |\stackrel{(\ref{L4H1})}{\leq} \left( \frac{n}{m_2}+ \sqrt{\beta}n \right)-
\left(\frac{\xi m_1}{6k'm_2}\right) \left(\frac{n}{4\Delta m_1}\right) \leq \frac{n}{m_2} ,$$
we have that $X'_{i, j} \in [0, 1]$, for each $1 \leq j \leq m_3$. Let
\begin{equation} \label{L4H3}
\mu = \sum_{j = 1}^{m_3} \frac{(1 + \xi/20)Y_j m_2}{k'n},
\end{equation}
and note that 
\begin{align} \label{L4H4}
\sum_{j = 1}^{m_3} Y_j &\leq \sum_{j = 1}^{m_3}|A \cap L_j| = \frac{1}{2}\left(\sum_{j = 1}^{m_3} |L_j| + s_j\right) \nonumber \\
&\stackrel{(\ref{L4H1}), (\ref{L4H2})}{\leq} \left(\frac{m_3 n}{2m_2} + \frac{m_3 \sqrt{\beta} n}{2}\right) + \left(\frac{|A| - |B|}{4} + \frac{\xi n}{40}\right) \leq \frac{m_3 n}{2m_2} + \frac{|A| - |B|}{4} + \frac{\xi n}{20}.
\end{align}
Note also that $Y_j \geq (n/(4\Delta m_1)) \times (m_1/2m_2)  = n/(8 \Delta m_2)$ for each $1 \leq j \leq m_1$. Thus we have $\mu \geq m_3/(8\Delta k') \gg (\log k')/\xi^2$. We now apply Lemma \ref{Azuma} with $\delta := \xi/20$ to obtain
$$\mathbb{P}\left[\sum_{j = 1}^{m_3} X'_{i, j} > (1 + \xi/20) \mu\right] \leq e^{-\frac{\xi^2 \mu}{1200}} \leq \frac{1}{3k'}.$$
It follows that with probability at least $1 - 1/3k'$,
\begin{align*}
X_i \stackrel{(\ref{L4H3})}{\leq} \frac{(1+ \xi/20)^2 }{k' } \sum_{i = 1}^{m_3} Y_j &\stackrel{(\ref{L4H4})}{\leq} \frac{(1 + \xi/20)^2}{k'}\left(\frac{m_3 n}{2m_2} + \frac{|A| - |B|}{4} + \frac{\xi n}{20}\right) \\
&\leq \frac{1}{2k'}\left(\frac{m_3 n}{m_2} + \frac{|A| - |B|}{2}\right) + \frac{\xi n}{6k'},
\end{align*}
which proves Claim \ref{claimb}. \medskip

By a similar argument we have that if $X'_i$ is the number of vertices of $\bigcup_{j = 1}^{m_3} S_j$ assigned to $2i$, then
$$\mathbb{P}\left[X'_i > \frac{1}{2k'}\left(\frac{m_3 n}{m_2} + \frac{|B| - |A|}{2}\right) + \frac{\xi n}{6k'}\right] \leq \frac{1}{3k'}$$
for every $1 \leq i \leq k'$. Taken together with Claim \ref{claimb} this implies that with probability at least $1/3$,
\begin{align} \label{L4H5}
&X_i \leq \frac{1}{2k'}\left(\frac{m_3 n}{m_2} + \frac{|A| - |B|}{2}\right) + \frac{\xi n}{6k'} \\
\text{ and } &X'_i \leq \frac{1}{2k'}\left(\frac{m_3 n}{m_2} + \frac{|B| - |A|}{2}\right) + \frac{\xi n}{6k'} \nonumber
\end{align}
for every $1 \leq i \leq k'$, and hence there exists an assignment such that (\ref{L4H5}) holds.

For each $m_3 < j \leq m_2$, let $D_j$ be the union of the \emph{first} $k_2$ small segments of $L_j$ and $S_j$ the union of the remaining small segments. We now assign the vertices of $\bigcup_{j = m_3 + 1}^{m_2} L_j$ using an algorithm similar to that for $\bigcup_{j = 1}^{m_3} L_j$, but in reverse order. That is, we first choose $1 \leq i''_0 \leq k'$ randomly and assign the vertices of $S_{m_2}$ using the sober algorithm, but with the roles of $A_i$ and $B_i$ exchanged for each $i$. Thus we assign the vertices of $B_{m_1}$ to $2i''_0 -1$, the vertices of $A_{m_1}$ to $2i''_0$, etc. Similarly we use the drunken algorithm to assign the vertices of $D_{m_2}$ (again with the roles of $A_i$ and $B_i$ exchanged for each $i$), and so on until we have assigned all the vertices up to $S_{m_3 + 1}$. 
(As before,  the initial vertex of any application of an algorithm is the successor of the final vertex of the previous application
of an algorithm.)
Finally we use the $2i'_2$-seeking algorithm to assign the vertices of the last drunken segment $D_{m_3 + 1}$. 
(Recall that the final vertex of the $2i'_2$-seeking algorithm is always $2i'_2$.)
Let $\overline{X_i}$ be the number of vertices of $\bigcup_{j = m_3 + 1}^{m_2} S_j$ assigned to to $2i-1$ and $\overline{X'_i}$ the number assigned to $2i$. By using a proof analogous to that of Claim \ref{claimb}, we can ensure that 
\begin{align} \label{L4H6}
&\overline{X_i} \leq \frac{1}{2k'}\left(\frac{(m_2 - m_3) n}{m_2} + \frac{|B| - |A|}{2}\right) + \frac{\xi n}{6k'} \\
\text{ and } & \overline{X'_i} \leq \frac{1}{2k'}\left(\frac{(m_2 - m_3) n}{m_2} + \frac{|A| - |B|}{2}\right) + \frac{\xi n}{6k'} \nonumber
\end{align}
for each $1 \leq i \leq k'$.

Note that 
$$|\bigcup_{j = 1}^{m_2} D_j| \leq m_2 \times \frac{\xi m_1}{6 k' m_2}\times (n/m_1+\beta n)
\leq \xi n/3k',$$ and hence in total we assign at most
\begin{align*}
&X_i + X'_i + |\bigcup_{j = 1}^{m_2} D_j| \\
\stackrel{(\ref{L4H5}), (\ref{L4H6})}{\leq} &\frac{1}{2k'}\left(\frac{m_3 n}{m_2} + \frac{|A| - |B|}{2}\right) + \frac{1}{2k'}\left(\frac{(m_2 - m_3) n}{m_2} + \frac{|B| - |A|}{2}\right) + \frac{2\xi n}{3k'} \\
= &\frac{1}{2k'}(1 + 4\xi /3)n
\end{align*}
vertices to $2i-1$ and at most $(1 + 4\xi /3)n/2k'$ vertices to $2i$, for each $1 \leq i \leq k'$. 
\medskip

\noindent
{\bf Completing the proof.}
This completes our definition of $f_2$. We now check that $f_2$ is a homomorphism. By properties (c) and (d) it suffices to show for each $i$ that whenever $A_i$ and $B_i$ (or $B_i$ and $A_{i+1}$) are assigned to vertices $1 \leq i', j' \leq 2k'$ of $C'$, then $i'j'$ is an edge of $C' \cup \{c'\}$.  Observe first that the sober, drunken and seeking algorithms all assign vertices in such a way that $i'j'$ is an edge of $C'$. Further, recall that the initial vertex of any application of an algorithm is the successor of the final vertex of the previous application of an algorithm. So if, for example, the vertices of $B_i$ are assigned to the final vertex $2i^*$ where 
$B_i$ is the final segment assigned in an application one of the algorithms, then the vertices of $A_{i+1}$
will be assigned to the initial vertex $2i^*+1$ in the next application of an algorithm.

The only pair this argument does not deal with is the pair $(B_j, A_{j+1})$, where $B_j$ is the last small segment of $D_{m_3}$ (and hence $A_{j+1}$ is the first small segment of $D_{m_3 + 1}$, and therefore the last to be assigned). Now by the definition of the seeking algorithm, the vertices of $B_j$ are assigned to the vertex $2i'_1$ of $C'$ and the vertices of $A_{j+1}$ are assigned to the vertex $2i'_2 $ of $C'$. Recalling that $c' = \{2i'_1, 2i'_2 \}$ we have that $f_2$ is indeed a homomorphism. 

Now consider $f = f_1 \circ f_2$. Since $f_1$ and $f_2$ are both homomorphisms we have that $f$ is a homomorphism. Property (b) implies that condition ($\beta_1$) holds.
By (c), every edge $xy$ not in $H[S]$ lies in a pair $(A_i,B_i)$ for some $i$. Thus, by definition of our three algorithms, $xy$
is mapped to an edge $\{2j-1, 2j\}$ by $f_2$ for some $j$. By definition of $f_1$, $\{2j-1, 2j \}$ is mapped to
$\{2j'-1,2j'\}$ by $f_1$ for some $j'$. Therefore, $f$ satisfies ($\beta _3$).

 To see that condition ($\beta_2$) also holds, recall that $f_1$ assigns to each vertex $2i-1$ of $C$ (and also to $2i$) exactly $\lceil (n_{2i-1} + n_{2i}) k_1/n \rceil$ vertices of $C'$. Hence $f$ assigns at most
\begin{align*}
\frac{1}{2k'}(1 + 4\xi /3)n \lceil (n_{2i-1} + n_{2i}) k_1/n \rceil &\leq (1 + 4\xi /3)(n_{2i - 1} + \xi n/10) + 
\frac{1}{2k'}(1 + 4\xi /3)n \\
&\leq n_{2i - 1} + \xi n
\end{align*}
vertices of $H$ to the vertex $2i-1$ of $C$, for each $1 \leq i \leq k$. Similarly $f$ assigns at most $n_{2i} + \xi n$ vertices of $H$ to the vertex $2i$. \endproof

\section{Completing the Proof} \label{SecCompleting}

In this section we use Lemmas \ref{LemmaForG} and \ref{LemmaForH} to prove Theorem \ref{mainthm}. We use the following 
definition and lemma from \cite{botphd}; these allow us to prove that $H$ embeds into $G$ by checking some relatively simple conditions. 

\begin{defin} \label{EpsilonCompat} Let $H$ be a graph on $n$ vertices, let $R$ be a graph on $[k]$, and let $R' \subseteq R$. We say that a vertex partition $V(H) = (W_i)_{i \in [k]}$ of $H$ is \emph{$\varepsilon$-compatible} with an integer partition $(n_i)_{i \in [k]}$ of $n$ and $R' \subseteq R$ if the following holds. For $i \in [k]$ let $S_i$ be the set of vertices in $W_i$ with neighbours in some $W_j$ with $ij \notin E(R')$ and $i \neq j$. Set $S := \bigcup_{i \in [k]} S_i$ and $T_i := N_H(S) \cap (W_i \backslash S)$. Then for all $i, j \in [k]$ we have that
\begin{itemize}
\item [\rm ($\gamma_1$)] $|W_i| = n_i$;
\item [\rm ($\gamma_2$)] $xy \in E(H)$ for $x \in W_i$ and $y \in W_j$ implies that $ij \in E(R)$;
\item [\rm ($\gamma_3$)] $|S_i| \leq \varepsilon n_i$ and $|T_i| \leq \varepsilon \cdot \min\{n_j \mid i$ and $j$ are in the same component of $R' \}$.
\end{itemize}
The partition $V(H) = (W_i)_{i \in [k]}$ is \emph{$\varepsilon$-compatible} with a partition $V(G) = (V_i)_{i \in [k]}$ of a graph $G$ and $R' \subseteq R$ if $V(H) = (W_i)_{i \in [k]}$ is \emph{$\varepsilon$-compatible} with $(|V_i|)_{i \in [k]}$ and $R' \subseteq R$.
\end{defin}

\begin{lemma}[\cite{botphd}, Lemma 3.12] \label{EmbedLemma} For all $d, \Delta, r > 0$ there is a constant $\varepsilon =  \varepsilon(d, \Delta, r)$ such that the following holds.
Let $G$ be a graph on $n$ vertices and suppose that $(V_i)_{i \in [k]}$ is a partition of $V(G)$. Suppose $R$ is an
 $(\varepsilon, d)$-reduced graph of $G$ on $V_1, \dots , V_k$ and that $R'$ is a subgraph of $R$ 
 whose connected components have size at most $r$. Assume that $(V_i,V_j)_G$ is an $(\eps, d)$-super-regular pair
 for every edge $V_iV_j \in E(R')$.
Further, let $H$ be a graph on $n$ vertices with maximum degree $\Delta(H) \leq \Delta$ that has a vertex partition $V(H) = (W_i)_{i \in k}$ which is $\varepsilon$-compatible with $V(G) = (V_i)_{i \in [k]}$ and $R' \subseteq R$. Then $H \subseteq G$. \endproof
\end{lemma}

Lemma \ref{EmbedLemma} is a consequence of the Blow-up lemma of  Koml\'os,  S\'ark\"ozy and Szemer\'edi~\cite{kss2}. 
We now prove Theorem \ref{mainthm}.

\proof[Proof of Theorem \ref{mainthm}] Firstly, note that it suffices to prove Theorem~\ref{mainthm} under the addition
assumption that $\eta \ll 1$.
We choose $\beta, n_0$ as well as additional constants $d, \varepsilon, \xi, \lambda$ as follows: Choose $d \ll \nu$ as required by Lemma \ref{LemmaForG} also ensuring that $d \ll 1/\Delta$. Then take $\varepsilon \leq \varepsilon(d, \Delta, 2)$ as in Lemma \ref{EmbedLemma}, ensuring also that $\varepsilon \ll d$. Finally choose
$$1/n_0 \ll \beta \ll \lambda \ll \xi \ll \varepsilon,$$
as required by Lemmas \ref{LemmaForG} and \ref{LemmaForH}.

Apply Lemma \ref{LemmaForG} to $G$ to obtain an integer $k$ such that $\xi \ll 1/k \ll \varepsilon$,   
a partition $(n_i)_{i = 1}^{2k}$ of $n$ and integers $1\leq i_1 \not =j_1 , i_2 \not =j_2\leq k$. 
Suppose $C$ is the cycle $12\dots (2k)1$ with the chord $c=\{2i_2, 2j_2\}$.
Next we apply Lemma \ref{LemmaForH} with the partition $(n_i)_{i = 1}^{2k}$ of $n$ as input to obtain a set $S \subseteq V(H)$ and a homomorphism $f: H \rightarrow C \cup \{c\}$, such that \begin{itemize}
\item[\rm (i)] $|S| \leq \xi n$;
\item[\rm (ii)] $|f^{-1}(i)| \leq n_i + \xi n$ for all $1 \leq i \leq 2k$;
\item[\rm (iii)] Every edge which is not in $H[S]$ is mapped to the edge $\{2i-1, 2i\}$, for some $1 \leq i \leq k$.
\end{itemize}
Let $W_i := f^{-1}(i)$ and $n'_i := |f^{-1}(i)|$ for each $i$, and note that $(n'_i)_{i = 1}^k$ is a partition of $n$. 
Condition (ii) together with Lemma~\ref{LemmaForG} imply that
there is a partition $A'_1, B'_1, A'_2, B'_2 , \dots ,A'_k, B'_k$ of $V(G)$ and a spanning subgraph $G'$ of $G$ which
satisfy conditions ($\alpha_1$)--($\alpha_5$).

Relabel these clusters $V_1, \dots , V_k$ such that $V_{2i-1}:=A'_i$ and $V_{2i}:=B'_i$ for all $1 \leq i \leq k$.
So $|V_i|=n'_i$ for all $1 \leq i \leq 2k$. Let $R$ be the $(\eps,d)$-reduced graph of $G'$ on
$V_1,\dots, V_k$ with the maximal number of edges. Hence ($\alpha _2$)--($\alpha _5$) imply that $R$ contains
the Hamilton cycle $C'=V_1V_2\dots V_{2k} V_1$ and the chord $c':=V_{2i_2}V_{2j_2}$ (we view $C'\cup\{c'\}$ as a copy
of $C\cup \{ c\} $ in $R$). Let $R'$ be the spanning subgraph of $R$ containing precisely the edges $V_{2i-1}V_{2i}$
for $1\leq i\leq k$. Note that ($\alpha _2 $) implies that $(V_{2i-1}V_{2i})_{G'}$ is an $(\eps,d)$-super-regular pair
for all $1\leq i\leq k$.

 We now check that the partition $V(H)=(W_i)_{i = 1}^{2k}$ is 
 $\varepsilon$-compatible with the partition $(V_i)_{i = 1}^{2k}$ and $R'\subseteq R$.
  We defined $(V_i)_{i = 1}^{2k}$ so that $|V_i| = |W_i|$ for each $1\leq i\leq 2k$ and hence condition ($\gamma_1$) of Definition \ref{EpsilonCompat} holds. Condition ($\gamma_2$) holds since $f: H \rightarrow C \cup \{c\}$ is a homomorphism and $C \cup \{c\}$ is a subgraph of $R$. Note that for all $1\leq i\leq 2k$, Lemma~\ref{LemmaForG} implies that 
  $$\eps n'_i \geq \eps( n_i -2k\xi n ) \geq \eps (n/3k-2k\xi n) \geq \eps n/4k \gg\xi n \stackrel{(i)}{\geq} |S|.$$
Furthermore, $|N_H (S) \cap (W_i \backslash S)|\leq \Delta |S|\leq \Delta \xi n \ll \eps n/4k \leq \eps n'_j$
for all $1\leq i,j\leq 2k$.
Thus,
  condition ($\gamma_3$)  holds.
Hence, Lemma \ref{EmbedLemma} implies that   $G'$ (and therefore $G$) contains $H$, as desired. \endproof

\section*{Acknowledgements}
We would like to thank J\'ozsef Balogh for helpful discussions, and Daniela K\"uhn and Deryk Osthus
for their comments on the manuscript. We also thank the referee for helpful comments and suggestions.

{\footnotesize \obeylines \parindent=0pt

\begin{tabular}{lll}

Fiachra Knox                        &\ &  Andrew Treglown \\
School of Mathematics						    &\ &  Faculty of Mathematics and Physics \\
University of Birmingham   					&\ &   Charles University\\
Birmingham                          &\ &  Malostransk\'e N\'am\v{e}st\'i 25\\
B15 2TT															&\ &  188 00 Prague \\
UK																	&\ &  Czech Republic\\
\end{tabular}
}

{\footnotesize \parindent=0pt

\it{E-mail addresses}:
\tt{knoxf@maths.bham.ac.uk}, \tt{treglown@kam.mff.cuni.cz}}


\begin{thebibliography}{99}
\bibitem{alon} N. Alon and R. Yuster, H-factors in dense graphs, \emph{J. Combin. Theory B}~{\bf 66} (1996), 269–-282.
\bibitem{bkt} J. Balogh, A.V. Kostochka and A. Treglown, On perfect packings in dense graphs, submitted.
\bibitem{botphd} J. B\"ottcher, \emph{Embedding large graphs--The Bollob\'as--Koml\'os conjecture and beyond},
PhD thesis, Technische Universit\"at M\"unchen, 2009.
\bibitem{bom} J. B\"ottcher and S. M\"uller, Forcing spanning subgraphs via Ore type conditions,
\emph{Electronic Notes in Discrete Mathematics},~{\bf 34} (2009), 255--259.
\bibitem{BPTW10} J. B\"ottcher, K. Preussmann, A. Taraz and A. W\"urfl,
Bandwidth, expansion, treewidth, separators and universality for bounded-degree graphs, 
\emph{European J. Combin.}~{\bf 31} (2010), 1217--1227.
\bibitem{BST08} J. B\"ottcher, M. Schacht and A. Taraz, Spanning $3$-colourable subgraphs of small bandwidth in dense graphs,
\emph{J. Combin. Theory B}~{\bf 98} (2008), 85--94.
\bibitem{bot} J. B\"ottcher, M. Schacht and A. Taraz, Proof of the bandwidth conjecture of Bollob\'as and Koml\'os,
\emph{Math. Ann.}~{\bf 343} (2009), 175--205.
\bibitem{chau} P. Ch\^au, An Ore-type theorem on hamiltonian square cycles, \emph{Graphs Combin.}, to appear.
\bibitem{ch} V. Chv\'atal, On Hamilton's ideals, \emph{J. Combin. Theory~B}~{\bf 12} (1972), 163--168.

\bibitem{corradi} K. Corr\'adi and A. Hajnal, On the maximal number of independent circuits in a graph, \emph{Acta Math. Acad. Sci.
Hungar.}~{\bf 14} (1964), 423--439.
\bibitem{csaba} B. Csaba, On embedding well-separable graphs, \emph{Discrete Math.}~{\bf 308} (2008), 4322--4331.
\bibitem{dirac} G.A. Dirac, Some theorems on abstract graphs, \emph{Proc. London Math. Soc.}~{\bf 2} (1952), 69--81.
\bibitem{posa} P. Erd\H{o}s, Problem 9, in: M. Fieldler (Ed.), \emph{Theory of Graphs and its Applications}, Czech. Acad. Sci. Publ., Prague, 1964, p. 159.
\bibitem{hs} A. Hajnal and E.~Szemer\'{e}di, Proof of a conjecture of Erd\H{o}s,
\emph{Combinatorial Theory and its Applications vol. II}~{\bf 4} 
(1970), 601--623.
\bibitem{han} H. H\`an, \emph{Einbettungen bipartiter Graphen mit kleiner Bandbreite}, Master's thesis, Humboldt-Universit\"at zu Berlin, Institut f\"ur Informatik,  2006.
\bibitem{huang} H. Huang, C. Lee and B. Sudakov, Bandwidth theorem for random graphs, \emph{J. Combin. Theory B} 
{\bf 102} (2012), 14--37. 
\bibitem{kier} H.A. Kierstead and A.V. Kostochka, An Ore-type Theorem on Equitable
Coloring, \emph{J. Combin. Theory B}~{\bf{98}}  (2008), 226--234.
\bibitem{komlos} J. Koml\'os, The Blow-up Lemma, \emph{Combin. Probab. Comput.}~{\bf 8} (1999), 161--176.
\bibitem{kss2} J. Koml\'{o}s, G.N. S\'ark\"{ozy} and E. Szemer\'{e}di, Blow-up
Lemma, \emph{Combinatorica}~{\bf 17} (1997), 109--123.
\bibitem{kss} J. Koml\'os, G.N. S\'ark\"ozy and E. Szemer\'edi, Proof of the Seymour conjecture for large graphs, \emph{Annals of Combinatorics}~{\bf 2} (1998), 43--60.
\bibitem{kss3} J. Koml\'os, G.N. S\'ark\"ozy and E. Szemer\'edi, Proof of the Alon--Yuster conjecture, \emph{Discrete Math.}~{\bf 235}
(2001), 255–-269.
\bibitem{KMOprep} D. K\"uhn, R. Mycroft and D. Osthus, An approximate version of Sumner's universal tournament conjecture,
\emph{J. Combin. Theory B}~{\bf 101} (2011), 415--447.
\bibitem{kuhn} D. K\"{u}hn and D. Osthus, Critical chromatic number and the
complexity of perfect packings in graphs, \emph{17th ACM-SIAM Symposium on 
Discrete Algorithms} (SODA 2006), 851--859.
\bibitem{kuhn2} D. K\"{u}hn and D. Osthus, The minimum degree threshold for perfect
graph packings, \emph{Combinatorica}~{\bf 29} (2009), 65--107.
\bibitem{kotore} D. K\"{u}hn, D. Osthus and A. Treglown, An Ore-type theorem for perfect packings in graphs,
\emph{SIAM J. Disc. Math.}~{\bf 23} (2009), 1335--1355.
\bibitem{KOT10} D. K\"uhn, D. Osthus and A. Treglown, Hamiltonian degree sequences in digraphs, \emph{J. Combin. Theory B}~{\bf{100}} (2010),
367--380.
\bibitem{ore} O.~Ore, Note on Hamilton circuits,
\emph{Amer. Math. Monthly} {\bf 67} (1960), 55. 
\bibitem{seymour} P. Seymour, Problem section, in \emph{Combinatorics: Proceedings of the British Combinatorial
Conference 1973} (T.P. McDonough and V.C. Mavron eds.), 201--202, Cambridge
University Press, 1974.
\bibitem{SVprep} B. Sudakov and J. Vondrak, A randomized embedding algorithm for trees, \emph{Combinatorica}~{\bf 30} (2010), 445--470.
\bibitem{reglem} E.~Szemer\'edi, Regular partitions of graphs, 
\emph{Probl\'emes Combinatoires et Th\'eorie des Graphes Colloques Internationaux
CNRS} {\bf 260} (1978), 399--401.





\end{thebibliography}
\end{document}